\documentclass{article}
\usepackage{amsfonts,amssymb}

\newtheorem{thm}{Theorem}
\newtheorem{dfn}{Definition}[subsection]
\newtheorem{prp}[dfn]{Proposition}
\newtheorem{lmm}[dfn]{Lemma}

\makeatletter

\@addtoreset{equation}{section}
\makeatother

\newcommand{\lbl}[1]%
{\label{#1}%
}
\newenvironment{eqnlbl}[1]%
{%
\begin{equation} \label{#1}}%
{\end{equation}}

\newcommand{\cst}{{\rm cst\,}}
\newcommand{\N}{{\rm I\hspace{-0.4ex}N}}
\newcommand{\qed}{\hspace*{2ex} \hfill $\square$}
\newcommand{\R}{{\rm I\hspace{-0.45ex}R}}
\newcommand{\Z}{{\sf Z\hspace{-0.8ex}Z}}
\newcommand{\T}[1]{{\cal T}\!\left[#1\right]}
\newcommand{\Tau}{{\cal T}}

\newcommand{\bx}[1]{\left[#1\right]}
\newcommand{\cf}{\underline}
\newcommand{\cl}{\overline}
\newcommand{\ds}[1]{{#1}^\natural}
\newcommand{\Ham}{\mbox{\rm H}}
\newcommand{\pt}{\hat}
\newcommand{\RC}{{\mbox{RC}}}
\newcommand{\hsh}{{\mbox{h-sh}}}
\newcommand{\vsh}{{\mbox{v-sh}}}
\newcommand{\Sp}{{\mbox{Sp}}}

\title{Collision probability\\
for random trajectories
in two dimensions}
\author{A. Gaudilli\`ere
    \thanks{Universit\`a Roma 3
    - gaudilli@mat.uniroma3.it}
    }  
\date{March $22^{{\rm nd}}$, 2007}

\begin{document}
\maketitle

\begin{abstract}
We give a lower bound
for the non-collision probability
up to a long time $T$
in a system of $n$ independent
random walks with fixed obstacles
on $\Z^2$.
By `collision' we mean collision between
the random walks
as well as collision
with the fixed obstacles.
We give an analogous result
for Brownian particles on the plane.
We also explain how this result
can be used to describe
in terms
of ``quasi random walks"
a diluted gas evolving
under Kawasaki dynamics
or simple exclusion.
\end{abstract}

\section{Results, motivations and strategy}

\subsection{Main results}

Consider $n$ particles
performing independent simple random walks
in continuous time on $\Z^2$:
with each particle we associate
a clock which,
independently of the other clocks,
rings following a Poisson process
of intensity 1,
and each time 
a particle's clock rings
this particle jumps
to one of its four nearest neighbours
with uniform law.
Assume now that these particles
evolve in the midst of 
a finite number of fixed obstacles,
{\em rectangles on $\Z^2$},
i.e. of the form
$$
\left([a,b]\times[c,d]\right)\cap\Z^2
\mbox{ with $a$, $b$, $c$, $d$ in $\R$},
$$
and say that a {\em collision} occurs
when some particle
is nearest neighbour
of one of the rectangular obstacles
or one of the other particles.
In this paper we give under a few
hypotheses on the initial configuration
of the system,
a lower bound
to the {\em non-collision probability}
up to time $T$,
for large $T$ and uniformly
in the initial configuration.
Denoting, for any $p\geq 1$,
by $d_p$ the distance associated
to the $p$-norm
$$
\|\cdot\|_p:(x,y)\in\R^2
\mapsto
\left\{
\begin{array}{ll}
\left(|x|^p+|y|^p\right)^{1/p}
& \mbox{if $p<+\infty$}\\
\sup(|x|,|y|)
& \mbox{if $p=+\infty$}
\end{array}
\right.,
$$
by $|E|$ the cardinality
of any finite set $E$,
by $a\wedge b$ the minimum
between the two real numbers $a$ and $b$,
and, for any $\ds{A} \subset \Z^2$
(in this paper the upper-index $\ds{}$
will identify 
the objects referring
to the lattice $\Z^2$),
by $\ds{\partial} \ds{A}$ its external
border
$$
\ds{\partial} \ds{A}:=
\left\{
z\in\Z^2\setminus \ds{A} : 
\exists z'\in \ds{A},
d_1(z,z')=1
\right\},
$$
our result reads
\begin{thm}
\lbl{thm1}
There exists a constant $c_0 \in ]0,+\infty[$ such
that for any $n\geq 2$ and $p\geq 2$
the following holds.
\par
Let $\ds{S}$ be a finite set of rectangles
$\ds{R}_1$, $\ds{R}_2$, \dots, $\ds{R}_{|\ds{S}|}$
on $\Z^2$ such that
$$
\left\{
\begin{array}{l}
|\ds{\partial} \ds{R}_1| + |\ds{\partial} \ds{R}_2| + \cdots
+ |\ds{\partial} \ds{R}_{|\ds{S}|}| \leq p\\
\inf_{i\neq j} d_\infty (\ds{R}_i,\ds{R}_j) > 3
\end{array}
\right.,
$$
and let $\cf{z}=(z_1,z_2,\dots,z_n) \in (\Z^2)^n$
be such that
$$
\left\{
\begin{array}{l}
\inf_{i\neq j} d_1 (z_i,z_j) > 1\\
\inf_{i,j} d_\infty (z_i,\ds{R}_j) > 3
\end{array}
\right..
$$
Denoting by $\ds{P}_{\cf{z}}$ the law
of $n$ independent simple random walks
in continuous time
$\ds{Z}_1$, $\ds{Z}_2$, \dots, $\ds{Z}_n$
starting from
$z_1$, $z_2$, \dots, $z_n$
and defining
$$
\ds{\Tau}_c:=
\inf
\left\{ t \geq 0:
\inf_{i\neq j} d_1 (\ds{Z}_i(t),\ds{Z}_j(t))
\wedge
\inf_{i,j} d_1 (\ds{Z}_i(t),\ds{R}_j)
= 1 \right\},
$$
we have
$$
\forall T \geq T_0,\quad
\ds{P}_{\cf{z}}\left(\ds{\Tau}_c>T\right)
\geq
\frac{1}{(\ln T)^\nu}\;,
$$
with
$$
\left\{
\begin{array}{l}
\nu =c_0 n^4 p^2 \ln p\\
T_0 =\exp\{\nu^2\}
\end{array}
\right..
$$
\end{thm}

\noindent
{\bf Remark:}
Since the perimeter of a rectangle on $\Z^2$
is at least 4,
the case $p=2$ corresponds to 
an empty set $S$.
In that case
the role played by the $\ds{R}_i$'s
is completely irrelevant.
\\

We will also prove a continuous version
of Theorem \ref{thm1}.
Consider $n$ spherical particles
of diameter 1,
centered at $n$ independent
planar Brownian motions
evolving in the midst
of `rectangles on the plane',
i.e. sets of the form
$$
[a,b]\times[c,d]
\mbox{ with $a$, $b$, $c$, $d$ in $\R$},
$$
and say that a collision
occurs when one of these particles
is tangent to
one of the rectangles
or to one of the other particles.
Denote, for $R$ a rectangle
on the plane,
by $|\partial R|$ its perimeter,
and, for $z$ a point in $\R^2$,
by $\pt{z}$ the part
of the plane occupied
by a particle centered
at $z$, i.e. the closed
ball centered at $z$
of diameter 1
for the distance $d_2$.
Then, in its continuous version
Theorem \ref{thm1} reads
\begin{thm}
\lbl{thm2}
There exists a constant $c_0 \in ]0,+\infty[$ such
that for any $n\geq 2$ and $p\geq 2$
the following holds.
\par
If $S$ is a finite set of rectangles
$R_1$, $R_2$, \dots, $R_{|S|}$
on the plane such that
$$
\left\{
\begin{array}{l}
|S| \leq p/4\\
|\partial R_1| + |\partial R_2| + \cdots
+ |\partial R_{|S|}| \leq p\\
\inf_{i\neq j} d_\infty (R_i,R_j) \geq 3
\end{array}
\right.,
$$
if $\cf{z}=(z_1,z_2,\dots,z_n) \in (\R^2)^n$
is such that
$$
\left\{
\begin{array}{l}
\inf_{i\neq j} d_2 (\pt{z_i}, \pt{z_j}) \geq 1\\
\inf_{i,j} d_\infty (\pt{z_i},R_j) \geq 3
\end{array}
\right.,
$$
then, denoting by $P_{\cf{z}}$ the law
of $n$ independent planar Brownian motions
$Z_1$, $Z_2$,~\dots, $Z_n$
starting from
$z_1$, $z_2$, \dots, $z_n$
and defining
$$
\Tau_c:=
\inf
\left\{ t \geq 0:
\inf_{i\neq j} d_2 (Z_i(t),Z_j(t))
\wedge
2\inf_{i,j} d_2 (Z_i(t),R_j)
= 1 \right\},
$$
we have
$$
\forall T \geq T_0,\quad
P_{\cf{z}}\left(\Tau_c>T\right)
\geq
\frac{1}{(\ln T)^\nu}\;,
$$
with
$$
\left\{
\begin{array}{l}
\nu = c_0 n^4 p^2 \ln p\\
T_0 = \nu^2
\end{array}
\right..
$$
\end{thm}

Actually dealing with this continuous
case is easier because it allows
for strong potential 
theoretic and stochastic techniques.
That is why we will first prove Theorem \ref{thm2}.
Theorem \ref{thm1} will be obtained
afterwards using the strong coupling
between Brownian motions and discrete random walks
built by Koml\'os, Major
and Tusn\'ady  (\cite{KMT1}, \cite{KMT2}).

\subsection{Motivations}
\label{mtvs}

The non-collision probability
estimated from below
in Theorems \ref{thm1} and \ref{thm2}
is a well known quantity
{\em in the one-dimensional case.}
In 1959 Karlin and McGregor
gave in \cite{KM} a determinant formula
to compute this one-dimensional
non-collision probability
without fixed obstacles
(the one-dimensional version
of our case $p=2$).
Their computation was based on
a reflection argument
which was after extended
in many different situations
(see for example \cite{Bi}, \cite{GZ},
\cite{HW}, \cite{Gr}, \cite{KOR}.)
In \cite{Gr} Grabiner
gave, for Brownian motions,
the asymptotic one-dimensional
continuous non-collision probability
without fixed obstacles
up to $T$ for large $T$ as
\begin{eqnlbl}{as1d}
c(\cf{\eta})
\left(1/\sqrt{T}\right)
^{n(n-1)/2}
\end{eqnlbl}
where $c(\cf{\eta})$ is an explicit
function of the initial configuration
$\cf{\eta}$.
He gave also some analogous
results for the non-collision
probability with one fixed obstacle
(the one-dimensional version
of our case $p=4$).

By analogy one can think that
in our two-dimensional case 
the non-collision probability
goes, for large T and at least
in the case $p= 2$,
as $(1/\ln T)^{n(n-1)/2}$
since $1/\ln T$ --
instead of $1/\sqrt{T}$ --
is the order of the probability
of not coming back to the origin
up to T.
The asymptotic $(1/\ln T)^{n(n-1)/2}$
is also the estimate
one would obtain
by assuming that the collisions
between different pairs of particles
are independent events.
We will turn back
at the end of the paper
to the question
of the right asymptotics,
we just note by now
that our lower bound is quite far
from this asymptotic
that one could expect,
and very far from the precision
of the one-dimensional asymptotic
given in (\ref{as1d}).
But the reflection argument
used in dimension one
does not apply
to the two-dimensional case --
at least not in the same direct way --
and we had to use a different approach
to get this weaker estimate.
As far as I know this is the first result
on the two-dimensional non-collision
probability.

Furthermore the lower bound
of Theorem \ref{thm1} -- 
uniform in the initial configuration,
given by $S$ and $\cf{z}$,
and valid for any $T$
larger than an explicit
(up to the constant $c_0$) $T_0$ --
is sufficient to describe 
in terms of ``quasi random walks"
a very diluted
lattice gas of density
$\rho=e^{-\alpha}$ 
with $\alpha \gg 1$
evolving under the Kawasaki dynamics
(or, in the particular case
of an infinite temperature,
under simple exclusion) in a large box
$\Lambda_\alpha \subset \Z^2$
with periodic boundary conditions
and exponentially large volume
$|\Lambda_\alpha| = e^{\Theta\alpha}$
for some positive $\Theta$.
This will be the object
of another publication
\cite{dHGNOS1}
and we just give here
the heuristic of the problem.

The Kawasaki dynamics
describes a system of particles
which evolve with exclusion
and interaction.
It is represented
by a Markov chain
on the space of configurations
$$
{\cal X} :=
\left\{\eta \in \{0;1\}^{\Lambda_\alpha}
:\:
\sum_{x\in\Lambda_\alpha} \eta(x) = \rho|\Lambda_\alpha|
\right\}
,$$
where 1 stands for an occupied site
and 0 for an empty one.
To each particle is associated
a clock which rings following
a Poisson process of intensity~1,
and each time its clock rings,
the particle ``tries'' to jump
to one of its four nearest neighbour sites,
randomly chosen with uniform law.
If this site is occupied
by another particle,
the particle the clock of which rang
stays in the same site.
If this site is vacant
the particle performs the jump
with probability
$$
p:= 
\exp \left\{
  -\beta \left[
    \Ham(\eta^2)-\Ham(\eta^1)
  \right]_+
\right\}
$$
where $\beta\geq 0$ stands
for the inverse temperature,
$\eta^1$ for the configuration
in which the system was
when the clock rang,
$\eta^2$ for the configuration
in which the system will be
if the particle jumps effectively     
in the chosen nearest neighbour site,
and $\Ham$ for the Hamiltonian
of the system defined by
$$
\Ham(\eta) :=
\sum_{
  \begin{array}{c}
    {\scriptstyle \{x;y\}\subset\Lambda_\alpha}\\
    {\scriptstyle \|x-y\|_1=1}
  \end{array}
}
-U\eta(x)\eta(y)
\qquad
\mbox{for any $\eta$ in ${\cal X}$}
,$$
with $-U<0$ the binding energy felt
by two neighbouring particles.
In particular a given particle
moves, under the Kawasaki dynamics,
like a simple random walk
as long as it is ``free'',
i.e., as long as its four nearest
neigbour sites are unoccupied.
When it collides with a cluster
or with another free particle
to form a new cluster,
it will be ``stuck'' for a while
inside this cluster,
but will be eventually ``re-emitted''.

Say now that the gas of the system
consists of all the free particles
and all the clusters smaller than
a finite given volume,
independent of $\alpha$.
Then the particles of the gas
will be re-emitted
after each collision at bounded
distance from the point where they
were clusterized.
And, roughly speaking,
they will perform
simple random walks perturbed
by this collision/re-emission process.
If the frequency of the collisions is ``low'',
then individual particles will perform
``quasi random walks" as introduced
in a simpler context in \cite{dHOS}.
By ``low frequency'' we mean
that the number of collisions
is non-exponentially
large in $\alpha$ in any time interval
of length $1/\rho=e^\alpha$.

Assume now that the system starts
from a measure $\mu$,
{\em which can be different
from the equilibrium measure,}
such that
one can {\em a priori} exclude
for very long times
(say exponentially large in $\alpha$)
any anomalous concentration
of the gas in any box of volume
$1/\rho$.
It means that with a very high probability
(say super-exponentially close to 1
in $\alpha$)
there will be no more
than $\lambda(\alpha)$
particles of the gas in any
of these boxes, with $\lambda$
a function which grows slowly
to $+\infty$ (for example
$\lambda(\alpha)=\ln\ln\alpha$):
for $c>0$ and with $G$ the subset
of $\Lambda_\alpha$ occupied
by the gas
$$
\limsup_{\beta\rightarrow +\infty}
\frac{1}{\beta}
\ln P_\mu\left(
  \exists t<e^{c\cdot\alpha},
  \exists \Lambda\subset\Lambda_\alpha,
  |\Lambda| \leq 1/\rho,
  |\Lambda \cap G| > \lambda(\alpha)
\right)
= -\infty.
$$
      
Then, on time scales $T= 1/\rho = e^{\alpha}$,
the clouds of potentially interacting particles,
i.e., at diffusive distance of order
$T^{1/2}=e^{\alpha/2}$
from each other, are ``virtually finite",
i.e., contain at most a number $n$ of particles 
of order $\lambda(\alpha)$.
Taking as first approximation
that the clusters are fixed obstacles,
an application of Theorem \ref{thm1}
and the strong Markov property
will give that the number of collisions
inside each cloud is higher than $\ln^{2\nu}T$
with a probability smaller than
$$
\left( 1 - \frac{1}{\ln^\nu T}\right)
^{\ln^{2\nu} T}
$$
which is super-exponentially small in $\alpha$.
So, with a probability super-exponential\-ly
close to 1, the number of collisions
inside a cloud of potentially interacting
particles is non-exponentially large
($\nu$ is linked to $n$ virtually finite),
this ensures in particular that
particles do not exhibit any super-diffusive
behaviour on time-scales $T$
(with a probability
super-exponentially close to 1)
and  
that the different clouds
do not interact with each other.
This allows us to conclude that the frequency
of the collisions is very low
and that individual particles
perform ``quasi random walks".

This is particularly relevant for the study
of the metastable regime
of the 
Kawasaki dynamics at low temperature ($\beta\gg 1$),
where $\alpha=\Delta\beta$,
with $\Delta>0$ an activity parameter.
The first paper \cite{dHOS}
dealing with this issue
introduced a simplified model
based on the assumption
that the interaction
between a given cluster
and its surrounding gas
``was like" an interaction
with a gas of independent random walks.
This was the basic assumption
which justified the introduction
of the so-called local Kawasaki dynamics
further studied in
\cite{dHNOS}, \cite{GOS}, \cite{NOS}, \cite{BdHN}. 
The description of the gas
in terms of ``quasi random walks''
is one of the key elements
which allow to extend the results
for simplified models to full Kawasaki dynamics
(see \cite{dHGNOS2}.)
And this was the original motivation of this paper.

\subsection{Strategy and notation}
\lbl{stratnnot}

Since we want to give lower bounds
which decrease slowly in $T$,
and since the probability
that a random walk
or a Brownian motion 
take less than time $T$
to travel a distance
of order $T^{\epsilon +1/2}$
decreases more than exponentially fast,
we will estimate
the probability of 
travelling such distances
without collision
to estimate the non-collision probability.
In practice we will do so
with $\epsilon = 1/2$:
other choices
would only affect
the value of the constant
$c_0$ appearing in
Theorems~\ref{thm1} and~\ref{thm2}.

So, in section \ref{c1obst}
we will estimate this probability
of travelling the distance $T$ 
without collision
for the continuous version of the system
in the simpler case
of a single particular fixed obstacle,
namely a spherical particle
fixed at the origin.
This is the crucial point of the proof
of our results
and uses basic potential theory on $\R^{2n}$.

In section \ref{1tomany}
we will prove Theorems \ref{thm1}
and \ref{thm2} in four steps.
We will first give some rough estimates
for the probability that particles
bypass the obstacles and increase
linearly the distance between
them without collision.
These estimates are simple
but somewhat technical,
they come more naturally
in the discrete case
and are easily adapted
to the continuous one.

Secondly, using these estimates
and some logarithmic
scale invariance property,
we will reinforce
the result of section \ref{c1obst}:
we will estimate, in the simpler case
of $n$ Brownian particles
and one fixed particle,
the probability of increasing
up to $T$
the distance between them,
without collision and avoiding
that any particle travel
a distance  $\alpha T$,
where $\alpha$ is a positive
constant depending only on $n$.

Thirdly, transferring the problem
on some `mesoscopic scale' $\sigma_0$
which lies between the `microscopic scale' 1
and the `macroscopic one' $T$
and is linked to the distance
between the fixed obstacles,
Theorem \ref{thm2} will then follow,
by induction on the number of obstacles,
from this reinforced result
and from the previous rough estimates.
We will eventually prove Theorem \ref{thm1}
using the strong coupling
between random walks
and Brownian motions
built by Koml\'os, Major
and Tusn\'ady in \cite{KMT1}
and \cite{KMT2}.
The use of this approximation
to deal with the discrete case
is at the origin of the different
expression for $T_0$
in Theorems \ref{thm1} and \ref{thm2}.

\subsubsection*{Notation}

In the whole paper we will use the
following notation.

In any dimension $d$
and for any $p\geq 1$,
we will denote by $d_p$ the distance
associated to the usual $p$-norm
on $\R^d$
and
by $B_p(\cf{z},r)$ --
where $\cf{z} \in \R^d$
and $r>0$ --
the open ball of center $\cf{z}$
and radius $r$
for the distance $d_p$.
The border and the closure
(for the topology
associated to these distances)
of any subset $A$ of $\R^d$
will be denoted $\partial A$
and $\cl{A}$.

For any $A \subset \R^2$
and any $r>0$
we define
$$
\ds{A} :=  A \cap \Z^2
$$
and
$$
\bx{A}_r := 
\bigcup_{z\in A}
\cl{B_\infty (z,r)}
.$$
If $A$ is restricted to a single point
$z \in \R^2$
we will write $\bx{z}_r$
instead of $\bx{\{z\}}_r$.

For $A \subset \R^2$
we also define
its {\em horizontal shadow}
$\hsh A$
and
its {\em vertical shadow}
$\vsh A$
as
\begin{eqnarray*}
\hsh A & := &
\left\{(x,y) \in \R^2 :\:
\exists\, (x', y') \in A,\,
x = x'\right\};\\ 
\vsh A & := &
\left\{(x,y) \in \R^2 :\:
\exists\, (x', y') \in A,\,
y = y'\right\}.
\end{eqnarray*}
The {\em circumscribed rectangle} of $A$,
denoted $\RC (A)$, 
is the intersection
of all the rectangles
$[a,b]\times[c,d]$ containing $A$.
If $A$ is a `rectangle on the plane', i.e., if $A=\RC(A)$, then we denote
by $|\partial A|$ its perimeter.
For $S$ any finite set of rectangles
on the plane,
we define
$$
\underline{S}
:=
\bigcup_{R\in S} R
.$$

For any $z$ in $\R^2$
we define
$$
\pt{z} := \cl{B_2\left(z,1/2 \right)}
$$
which is the region occupied
by a spherical particle
with unitary diameter
centered at $z$.
For the discrete case we define
the following analogue:
$$ 
\bx{z} := \bx{z}_\frac{1}{2}
.$$

For any
$\cf{z} = (z_1, z_2, \dots, z_n)$ in $(\R^2)^n$
and any
$S$ finite set of rectangles
$R_1$, $R_2$,~\dots, $R_{|S|}$
on the plane,
we define two measures
$w_S(\cf{z})$ and $\ds{w}_S (\cf{z})$
of the distances between the particles
centered at $z_1$, $z_2$,~\dots, $z_n$
and the rectangles,
one for the continuous case,
the other for the discrete one:
$$
w_S(\cf{z}) :=
\inf_{i\neq j} d_\infty \left(\pt{z_i},\pt{z_j}\right)
\wedge
\inf_{i,j} d_\infty \left(\pt{z_i},R_j\right)
$$
and
$$
\ds{w}_S (\cf{z}):=
\inf_{i\neq j}
d_\infty \left(\bx{z_i},\bx{z_j}\right)
\wedge
\inf_{i,j}
d_\infty \left(\bx{z_i},R_j\right)
.$$
When there will be no ambiguity
on the set $S$ which these quantities
are referred to,
we will omit the index $_S$.
We also define,
with $O$ the origin of the plane,
$$
\delta(\cf{z}) :=
\inf_{i\neq j} d_2 (z_i,z_j)
\wedge
\inf_{i} d_2(z_i,O).
$$
Note that in this last definition
we take into account the distances
between {\em the centers} of the particles,
and not between the particles themselves.

For $Z_1$, $Z_2$, \dots, $Z_n$
$n$ independent planar Brownian motions,
we will denote by $\cf{Z}$ the $2n$-dimensional
Brownian motion
$$
\cf{Z} :=
\left( Z_1, Z_2, \cdots, Z_n \right)
$$
and for any $t \geq 0$
we define a ``maximal individual elongation"
up to $t$:
$$
\rho(t) :=
\sup_i
\sup_{s\leq t}
d_2(Z_i(s),Z_i(0)).
$$
In the same way, if $\ds{Z}_1$, $\ds{Z}_2$,
\dots, $\ds{Z}_n$ are $n$ independent random walks
in continuous time,
we will denote by $\ds{\cf{Z}}$ the process
$$
\ds{\cf{Z}} :=
\left( \ds{Z}_1, \ds{Z}_2, \cdots, \ds{Z}_n \right)
.$$

The {\em first collision time} $\Tau_c$
(respectively $\ds{\Tau}_c$ in the discrete case)
is defined for a given set $S$ 
(respectively $\ds{S}$)
of rectangles
as in Theorem \ref{thm2}
(respectively Theorem~\ref{thm1}.)
When we want to stress
the dependence on $S$ (or $\ds{S}$),
we will write
$\Tau_{c;S}$
(respectively $\ds{\Tau}_{c;\ds{S}}$.)
For any $A \subset (\R^2)^n$
and $b\geq 0$
we define
the stopping times
\begin{eqnarray*}
\T{A} &
:= & 
\inf
\left\{ t \geq 0 :
\cf{Z}(t) \in A
\right\}
,\\
\T{\rho \geq b} &
:= &
\inf
\left\{ t \geq 0 :
\rho(t) \geq b
\right\},
\end{eqnarray*}
in the same way we define
$\T{\delta \geq b}$, $\T{w \geq b}$
and $\T{\ds{w} \geq b}$,
we define also,
with $O$ the origin
of the plane,
$$
\Tau_{c;\pt{O}} 
:=
\left\{ t \geq 0:
\inf_{i\neq j} d_2 (Z_i(t),Z_j(t))
\wedge
\inf_{i,j} d_2 (Z_i(t),O)
= 1 \right\}.
$$
This last stopping time
is the extension
of the first collision time
to a situation in which
the set of fixed obstacles
is not made of rectangles
but of a single fixed particle
$\pt{O}$, centered at the origin $O$.

We will denote by $P_{\cf{z}}$
(respectively $\ds{P}_{\cf{z}}$)
the law of $n$ independent
planar Brownian motions
(respectively $n$ continuous time
planar random walks)
starting from 
$\cf{z}$ in $(\R^2)^n$.

We will use the notation
$$
\R_+ := \left\{x\in \R :\:
x \geq 0 \right\}
,$$
the convention
$$
\inf \emptyset := +\infty
,$$
and in all our computations
`$\cst$' will denote
a positive constant independent
of any parameter,
and the value of which
can change from line to line.

\section{Brownian motions with a single obstacle}
\lbl{c1obst}

\subsection{The key lemma}

In this section we study the simpler
case of $n$ Brownian particles
and a single fixed obstacle $\pt{O}$.
The following lemma
is the key point
of the proof of Theorems~\ref{thm1}
and~\ref{thm2}.

\begin{lmm}
\lbl{klmm}
For any $n\geq 2$, any $a \geq 2$ and any $T>0$,
if $\cf{z} \in (\R^2)^n$
is such that $\delta(\cf{z}) \geq a$,
then
\begin{eqnlbl}{kest}
\begin{array}{l}
P_{\cf{z}} \left(
\Tau_{c;\pt{O}} >
\T{\partial B_2(\cf{z},T)}
\right)\\
\quad \geq \quad
\left[
\left(
\frac{\ln a}{\ln (a +T)}
\right)
^{n}
\left(
\frac{\ln a}{\ln (a +\sqrt{2}T)}
\right)
^{\frac{{\scriptstyle n(n-1)}}{{\scriptstyle 2}}}
\right]
^{\frac{{\scriptstyle 1}}
{{\scriptstyle 1 -
\cos \frac{{\scriptstyle \pi}}{{\scriptstyle 2n}}}}}
,\end{array}
\end{eqnlbl}
so that, in particular,
for any $\epsilon>0$ 
and $b\geq a^{1+\epsilon}$,
\begin{eqnlbl}{lnsclinv}
P_{\cf{z}} \left(
\Tau_{c;\pt{O}} >
\T{\partial B_2(\cf{z},b)}
\right)
\geq
\left(
\frac{\ln a}{\ln b}
\right)
^{c_\epsilon n^4}
,
\end{eqnlbl}
where $c_\epsilon$
is a constant
which depends only on $\epsilon$.
\end{lmm}

\noindent
{\bf Proof:}
We denote by $E$ the Euclidean space $(\R^2)^n$
dressed with the usual scalar product,
and introduce
the subspaces of codimension
two $F_k$,
corresponding to the superposition
of two particles:
the $F_k$'s are all the subspaces
of the first kind
$$
F_k :=
\left\{
\cf{z} = (x_1,y_1,x_2,y_2,\dots,x_n,y_n)
\in E :
x_i = x_j,
y_i = y_j
\right\}
$$
for some $1\leq i < j \leq n$,
or of the second kind
$$
F_k :=
\left\{
\cf{z} = (x_1,y_1,x_2,y_2,\dots,x_n,y_n)
\in E :
x_i = 0,
y_i = 0
\right\}
$$
for some $1 \leq i \leq n$
(this last case corresponding to the superposition
of the $i^{th}$ particle
and the fixed one.)
So, $k$ is an integer index
going from 1 to $m$ with
$$
m:= \frac{n(n-1)}{2} + n
.$$ 
Observe that, denoting by $r_k$
the Euclidean distance (in $E$) to $F_k$,
we have, for any $\cf{z} = (z_1, z_2, \dots, z_n)$
in $E$,
$$
r_k (\cf{z}) =
\left\{
\begin{array}{ll}
\frac{1}{\sqrt{2}} \times d_2(z_i,z_j)
& \mbox{if $F_k$ is 
associated with the indices $i$ and $j$
($1^{st}$ kind);}\\
1 \times d_2(z_i,O) 
& \mbox{if $F_k$ is 
associated with the single index $i$
($2^{nd}$ kind).}
\end{array}
\right.
$$  
Calling $\alpha_k$ the inverse
of the corresponding coefficient
(so that $\alpha_k \in \{\sqrt{2};1\}$)
and defining the subsets of $E$
\begin{eqnarray*}
A_k & := &
\left\{
\cf{z} \in E :
\alpha_k r_k (\cf{z}) \leq 1 
\right\}\\
A & := & \bigcup_{1\leq k \leq m} A_k
\end{eqnarray*}
we get $\Tau_{c;\pt{O}} = \T{A}$
if the process starts
from some $\cf{z} \not\in A$.
So, with
$$
B^0:= B_2(\cf{z}^0,T)
$$
for a given $\cf{z}^0$ 
such that $\delta(\cf{z}^0) \geq a$,
we have to estimate
$P_{\cf{z}^0}(\T{A} > \T{\partial B^0})$
from below.

The function of the starting point
$$
h(\cdot):= P_{{\textstyle \cdot}}
\left( \T{A} > \T{\partial B^0} \right)
$$
is harmonic on $B^0 \setminus A$
and satisfies
$$
h|_{\partial B^0 \setminus A}
\equiv 1
\qquad
\mbox{and}
\qquad
h|_{\partial A}
\equiv 0
.$$
Assume now that $g$
is a subharmonic and non-negative function
on $E\setminus A$
that can be continuously extended
to get
$$
g|_{\partial A} \equiv 0
.$$
Then one gets
$$
h
\geq
\frac{g}{\sup g(B^0)}
$$
and this gives in $\cf{z}^0$
\begin{eqnlbl}{gest}
P_{\cf{z}^0}(\T{A} > \T{\partial B^0})
\geq
\frac{g(\cf{z}^0)}{\sup g(B^0)}
\end{eqnlbl}
(One can get the same result
applying Doob's Theorem 
to a family of stopped processes
obtained from the local submartingale
$g(\cf{Z})$.)

Now we look for such a function $g$
to get (\ref{kest})
as a consequence of (\ref{gest}).
Let us try with
$$
g = \prod_{1\leq k \leq m} f_k
$$
where each $f_k$ is an
increasing $C^2$ function
of $r_k$ such that
$$
f_k|_{\partial A_k}
\equiv 0
.$$
Note that in such conditions
we have
$$
\nabla f_k =
\left\|\nabla f_k\right\|
\nabla r_k
$$
so that, on $E \setminus A$
\begin{eqnarray*}
\frac{\Delta g}{g}
& = &
\sum_k \frac{\Delta f_k}{f_k}
+ \sum_{k\neq l} \frac{\nabla f_k}{f_k}
\cdot
\frac{\nabla f_l}{f_l} \\
& = &
\sum_k \frac{\Delta f_k}{f_k}
+ \sum_{k,l} \left\|\frac{\nabla f_k}{f_k}\right\|
\left\|\frac{\nabla f_l}{f_l}\right\| 
\nabla r_k \cdot \nabla r_l
- \sum_k \left\|\frac{\nabla f_k}{f_k}\right\|^2
.\end{eqnarray*}
Denoting, 
for any $\cf{z}$ in $E \setminus A$,
by $W(\cf{z})$ the {\em non-negative coordinates}
vector
$$
W(\cf{z}) := 
\left(
\left\|
\frac{\nabla f_k}{f_k}
(\cf{z})
\right\|
\right)_{1 \leq k \leq m}
\in \R_+^m
,$$
and by $Q(\cf{z})$ the $m$-dimensional real
symmetric matrix
$$
Q(\cf{z}):=
\left(
\nabla r_k (\cf{z}) \cdot \nabla r_l (\cf{z})
\right)_{1\leq k,l \leq m}
\in {\cal M}_m (\R)
,$$
we get, on $E\setminus A$,
\begin{eqnlbl}{laplest}
\frac{\Delta g}{g} =
\sum_k \frac{\Delta f_k}{f_k}
+ \left\|W\right\|_2^2
\left(
\frac{^tW}{\|W\|_2}
Q
\frac{W}{\|W\|_2}
-1
\right) 
,\end{eqnlbl}
where $^tW$ stand for
the line matrix obtained
by transposition from $W$.

We call {\em collision correlation}
the last factor
in (\ref{laplest})
(if the index $k$ described
a subset of $\{1;\dots;m\}$ such
that the associated $F_k$ correspond
to independent collisions, like for example 
between the first and second particle 
and between the third and the fourth one,
then this factor would be equal to 0).
This collision correlation
can be estimated from below
by $\gamma -1$ with
$$
\gamma
:=
\inf_{\cf{z} \in E \setminus A}
\inf_{
\begin{array}{c}
{\scriptstyle V \in \R_+^m}\\
{\scriptstyle \|V\|_2 = 1}
\end{array}
}
{^tV}Q(\cf{z})V
.$$
We claim
\begin{lmm}
\lbl{gammest}
For any $n \geq 2$
$$
\gamma \geq 1 - \cos\frac{\pi}{2n}
.$$
\end{lmm}
We postpone the proof of this result
to the next subsection
and note that, since (\ref{laplest})
implies that
$$
\frac{\Delta g}{g} \geq
\sum_k \frac{\Delta f_k}{f_k}
+ (\gamma -1)
\left\|\frac{\nabla f_k}{f_k}\right\|^2
,$$
a sufficient condition to get
the subharmonicity of $g$
is that all the $f_k$'s are
solution of the differential inequality
on $E\setminus A$:
\begin{eqnlbl}{diffineqn}
\frac{\Delta f}{f}
+ (\gamma-1)
\left\|\frac{\nabla f}{f}\right\|^2 
\geq 0
.\end{eqnlbl}
Since Lemma \ref{gammest} states
that $\gamma >0$, it is straightforward
to check that a positive $f$ is solution
of (\ref{diffineqn}) if and only if
$f^\gamma$ is subharmonic.
This shows (recall that the $r_k$'s 
measure the distance to
subspaces of codimension 2)
that we can choose
for every $k$
$$
f_k = (\ln \alpha_k r_k)^{\frac{1}{\gamma}}
.$$
By (\ref{gest}) we get then
\begin{eqnarray*}
P_{\cf{z}^0}(\T{A} > \T{\partial B^0})
& \geq &
\frac{1}{\sup g(B^0)}
\prod_k \left(
\ln \alpha_k r_k (\cf{z}^0)
\right)^{\frac{1}{\gamma}}\\
& \geq &
\prod_k \left(
\frac{\ln \alpha_k r_k (\cf{z}^0)}
{\ln(\alpha_k r_k (\cf{z}^0) + \alpha_k T)}
\right)^{\frac{1}{\gamma}}
,\end{eqnarray*}
and, since for any $k$
$$
x \in [2,+\infty[
\longmapsto 
\frac{\ln x}{\ln( x+ \alpha_k T)}
$$
is an increasing function,
this,
with the estimate
of $\gamma$ given by
Lemma \ref{gammest},
concludes the proof.
\qed

\subsection{Estimating the collision correlation}

We prove now Lemma \ref{gammest}
and keep the same notation
as in the previous subsection.
Any $\nabla r_k (\cf{z})$
belongs to the orthogonal
of $F_k$,
which is of dimension 2,
and the direction of the $\nabla r_k$'s   
depends on the point $\cf{z}$
where there are computed.
As a consequence $Q(\cf{z})$
depends strongly on $\cf{z}$.
But, as a matter of fact,
$\gamma$ can be estimated
from similar quantities
{\em computed for subspaces
of codimension 1.}

\subsubsection{Reducing the codimension}

To show this property we introduce some more notation.
Let us denote by
$$
(e^x_1, e^y_1, e^x_2, e^y_2, \dots, e^x_n, e^y_n)
$$
the canonical base of $E = (\R^2)^n$.
For $*$ any of
the two letters $x$ and $y$
we define
$$
F_k^* :=
\left\{
\cf{z} = (x_1,y_1,x_2,y_2,\dots,x_n,y_n)
\in E :
*_i = *_j
\right\}
$$
if $F_k$ corresponds to the superposition
of the $i^{th}$ and $j^{th}$
particles and
$$
F_k^* :=
\left\{
\cf{z} = (x_1,y_1,x_2,y_2,\dots,x_n,y_n)
\in E :
*_i = 0
\right\}
$$
if $F_k$ corresponds to the superposition
of the $i^{th}$ and the fixed one.
Calling $p_k^*$ the orthogonal projection
on $F_k^*$ we set
$$
u_k^* (\cf{z}) :=
\left\{
\begin{array}{ll}
\frac{\cf{z} - p_k^* \cf{z}}
{\left\|\cf{z} - p_k^* \cf{z}\right\|}
& \mbox{if $\cf{z} \in E\setminus F_k^*$}\\
0
& \mbox{if $\cf{z} \in F_k^*$}
\end{array}
\right.
.$$
Note that
$\|u_k^*(\cf{z})\| \in \{0;1\}$
and
$u_k^*(\cf{z})$ is collinear to some
\begin{eqnlbl}{dfnf}
f^*_{i,j}:=
\frac{e_j^* - e_i^*}{\sqrt{2}}
\end{eqnlbl}
if $F_k$ is of the first kind
or collinear to some
$e^*_i$ if $F_k$ is of the second kind.
It is also straightforward to check that
\begin{eqnlbl}{ortho}
\forall k,l \in \{1;\dots;m\},\:
u_k^x \cdot u_l^y \equiv 0.
\end{eqnlbl}
and 
\begin{eqnlbl}{dfnlambda}
\begin{array}{l}
\forall k \in \{1; \dots; m\},\:
\forall \cf{z} \in E\setminus A,\:
\exists ! (\lambda^x_k,\lambda^y_k)
\in \left[0, \|u_k^x(\cf{z})\| \right]
\times
\left[0, \|u_k^y(\cf{z})\| \right],\\
\quad
(\lambda^x_k)^2 + (\lambda^y_k)^2 = 1
\quad and \quad
\nabla r_k (\cf{z}) = \lambda^x_k u_k^x(\cf{z})
+ \lambda^y_k u_k^y(\cf{z}).
\end{array}
\end{eqnlbl}
We also define,
for any $\cf{z} \in E$,
$$
\gamma^* (\cf{z})
:=
\inf
\left\{
{^tV^*}Q^*(\cf{z})V^* \in \R :\:
V^* \in \prod_k \left[0, \|u_k^*(\cf{z})\| \right]
\subset \R_+^m,\:
\|V^*\|_2 = 1
\right\}
,$$
with
$$
Q^*(\cf{z}):=
\left(
u_k^* (\cf{z}) \cdot u_l^* (\cf{z})
\right)_{1\leq k,l \leq m}
\in {\cal M}_m (\R)
.$$

Now for $\cf{z}$ in $E\setminus A$
and
$$
V = (\mu_1, \mu_2, \dots, \mu_m) \in \R_+^m
$$
such that $\|V\|_2 = 1$,
writing
$$
V^* (\cf{z}) := (\mu_1\lambda^*_1,\dots,\mu_m\lambda^*_m)
\in \prod_k \left[0, \|u_k^*(\cf{z})\| \right]
,$$
where the $\lambda^*_k$'s
are defined by (\ref{dfnlambda}),
we have, using (\ref{ortho}),
\begin{eqnarray*}
{^tV} Q V
& = & \left\| \mu_1 \nabla r_1 + \cdots + \mu_m \nabla r_m \right\|^2\\
& = & \left\| \sum_k \mu_k\lambda^x_k u_k^x
      + \sum_k \mu_k\lambda^y_k u_k^y \right\|^2\\
& = & \left\| \sum_k \mu_k\lambda^x_k u_k^x \right\|^2
      + \left\| \sum_k \mu_k\lambda^y_k u_k^y \right\|^2\\
& = & {^tV^x} Q^x V^x + {^tV^y} Q^y V^y\\
& \geq & \left\|V^x\right\|^2_2 \gamma^x 
+ \left\|V^y\right\|^2_2 \gamma^y 
.\end{eqnarray*}
By the symmetry of the definitions
$$
\gamma_1 := \inf_{\cf{z} \in E} \gamma^x (\cf{z})
= \inf_{\cf{z} \in E} \gamma^y (\cf{z})
.$$
The equations in (\ref{dfnlambda})
give also
$$
\left\|V^x\right\|^2_2 + \left\|V^y\right\|^2_2
\equiv \left\|V\right\|^2_2 =1
$$
and we can conclude
${^tV} Q V \geq \gamma_1,$
so that
$\gamma \geq \gamma_1$
and we just have
to give a lower bound to $\gamma_1$,
i.e., a uniform lower bound
to $\gamma^x$ or $\gamma^y$, say $\gamma^x$.

\subsubsection{Estimating $\gamma^x$}

Any $u^x_k(\cf{z})$
which appears in the definition
of $\gamma^x(\cf{z})$
depends only on ``the side of the hyperplane $F^x_k$
where $\cf{z}$ lies."
As a consequence the function $\gamma^x$
is constant on any connected component
of $E\setminus \bigcup_k F^x_k$.
It is easy to see that
for any $\cf{z}' \in \bigcup_k F^x_k$
the infimum which defines
$\gamma^x(\cf{z}')$ is computed
on a set contained
in the one used
to compute $\gamma^x(\cf{z})$
for some $\cf{z}$ in $E \setminus \bigcup_k F^x_k$. 
So, to give a lower bound
to $\gamma^x(\cf{z})$
uniform in $\cf{z}$
we can assume that 
$$\cf{z} = (x_1, y_1, \dots, x_n, y_n)
\in
E\setminus \bigcup_k F^x_k
,$$
and this means that
0 and the coordinates
$x_1$, \dots, $x_n$ are $n+1$
distinct numbers.
Without loss of generality
we can then assume that
\begin{eqnlbl}{ord}
x_1 < x_2 < \cdots < x_n
\end{eqnlbl}   
and we have to show
$$
\gamma^x(\cf{z}) 
\;=\;
\inf_{
\begin{array}{c}
{\scriptstyle (\mu_1, \dots, \mu_m) \in \R^m_+}\\
{\scriptstyle \mu_1^2+ \cdots + \mu_m^2 =1}
\end{array}
}
\left\| \mu_1 u^x_1(\cf{z})
+ \cdots
+ \mu_m u^x_m (\cf{z}) \right\|^2
\;\geq\;
{1-\cos\frac{\pi}{2n}}
\;.$$

We will prove this lower bound in two steps.
First we will show that we can extract
from the family of the $m$ vectors $u^x_k(\cf{z})$
a family of $n$ vectors $v_1$, \dots, $v_n$
such that
\begin{eqnlbl}{nvect}
\gamma^x(\cf{z}) 
= \inf_{
\begin{array}{c}
{\scriptstyle (\mu_1, \dots, \mu_n) \in \R^n_+}\\
{\scriptstyle \mu_1^2+ \cdots + \mu_n^2 =1}
\end{array}
}
\left\| \mu_1 v_1
+ \cdots
+ \mu_n v_n \right\|^2
.\end{eqnlbl} 
Secondly we will show that this infimum
is greater than or equal to
$$
\inf \Sp (Q_n) = 1-\cos\frac{\pi}{2n}
\;,$$
the smallest eigenvalue of $Q_n$,
defined by
\begin{eqnlbl}{dfnQn}
Q_n :=
\left(
   \begin{array}{ccccc}
      1  & -\frac{1}{\sqrt{2}} & 0 & \cdots & 0 \\
      -\frac{1}{\sqrt{2}} & 1 & -\frac{1}{2} & \ddots & \vdots\\
      0 & -\frac{1}{2} & 1 & \ddots & 0 \\
      \vdots& \ddots & \ddots & \ddots & -\frac{1}{2} \\
      0 & \cdots & 0 & -\frac{1}{2} & 1 \\
   \end{array}
\right)
\in
{\cal M}_n(\R)
.\end{eqnlbl}

\subsubsection{From $m$ to $n$ vectors}

Defining
$$
q:= (n+1) \wedge \inf \left\{i \in \N :\: x_i >0 \right\} 
$$
we have (recall (\ref{ord}) and (\ref{dfnf}))
\begin{eqnarray*}
&&\left\{ u_k^x (\cf{z}) :\:
1\leq k \leq m \right\}\\
&& \quad = \: \left\{ f_{i,j}^x :\: 1\leq i < j \leq n \right\}
\;\cup\;
\left\{-e^x_i :\: 1 \leq i<q \right\}
\;\cup\;
\left\{e^x_i :\: q \leq i \leq n \right\} 
.\end{eqnarray*}
We define for $1 \leq i \leq n$
$$
v_i:=
\left\{
   \begin{array}{ll}
      f^x_{i,i+1} & \mbox{if $1 \leq i \leq q-2$}\\
      - e^x_{i}   & \mbox{if $i = q-1$}\\
      e^x_i       & \mbox{if $i = q$}\\
      f^x_{i-1,i} & \mbox{if $q+1 \leq i \leq n$}
   \end{array}
\right.
$$
It is easy to see that, for any $v_i$ and $v_j$,
with $i \leq j$,
\begin{eqnlbl}{scalp}
v_i \cdot v_j =
\left\{
   \begin{array}{rl}
      0 & \mbox{if $j  \geq i+2$}\\
      -\frac{1}{2} & \mbox{if $j = i+1 < q-1$}\\
      -\frac{1}{\sqrt{2}} & \mbox{if $j = i+1 = q-1$}\\
      0 & \mbox{if $j = i+1 = q$}\\
      -\frac{1}{\sqrt{2}} & \mbox{if $j = i+1 = q+1$}\\
      -\frac{1}{2} & \mbox{if $j = i+1 > q+1$}\\
      1 & \mbox{if $j = i$}
   \end{array}
\right.
\end{eqnlbl}
and any $u^x_k(\cf{z})$ is a {\em non-negative}
linear combination of the $v_i$'s:
$$
u_k^x(\cf{z}) =
\sum_i \lambda_{k,i} v_i
\quad \mbox{with} \quad
(\lambda_{k,i})_{1\leq i \leq n}
\in \R_+^n
.$$ 
Note that this equation
gives together with (\ref{scalp})
\begin{eqnlbl}{sumsqrs}
\sum_i \lambda_{k,i}^2
\;=\; 1 - \sum_{i\neq j} \lambda_{k,i}\lambda_{k,j} (v_i \cdot v_j)
\;\geq\; 1
.\end{eqnlbl}
Now for any
$$
(\mu_k)_{1\leq k \leq m} \in \R_+^m
$$
\mbox{such that}
$$
\sum_k \mu_k^2 = 1
$$
we have
$$
\left\|\sum_k \mu_k u_k^x(\cf{z})\right\|^2
\;=\; \left\|\sum_k \mu_k \sum_i \lambda_{k,i} v_i\right\|^2\\
\;=\; \left\|\sum_i \left(\sum_k \lambda_{k,i} \mu_k\right) v_i \right\|^2
.$$
Set
$$
s^2 := \sum_i \left(\sum_k \lambda_{k,i} \mu_k\right)^2
$$
and
$$
\left(\mu_i'\right)_{1\leq i\leq n} :=
\left(   
\frac{1}{s} \sum_k \lambda_{k,i} \mu_k
\right)_{1 \leq i \leq n}
\in \R_+^n
,$$
we get
$$
\sum_i \mu_i'^2 = 1
$$
and
$$
\left\|\sum_k \mu_k u_k^x(\cf{z})\right\|^2
= s^2 
\left\|\sum_i \mu_i' v_i \right\|^2
.$$
So, provided that $s^2 \geq 1$,
we get (\ref{nvect}).
But, using the fact that the $\lambda_{k,i}$'s
and the $\mu_k$'s 
are non-negative and
using (\ref{sumsqrs}), we have
$$
s^2
\;\geq\;
\sum_i \sum_k \lambda_{k,i}^2 \mu_k^2
\;=\;
\sum_k \mu_k^2 \sum_i \lambda_{k,i}^2
\;\geq\;
\sum_k \mu_k^2
\;=\;
1
$$
and this concludes our first step.

\subsubsection{The eigenvalues of $Q_n$}

Now it is easy to see that relations (\ref{scalp})
give that
$\gamma^x(\cf{z})$ is greater than or equal to
(equal to in the case $q=1$)
$$
\epsilon_n :=
\inf_{X \in K}
{^tX}Q_nX
$$
where $K$ is the closure of
$$
{\cal O} := 
\left\{ X = (x_1, \dots, x_n) \in ]0,+\infty[^n
:\:
\|X\|_2^2=1\right\}
$$
and $Q_n$ is defined in (\ref{dfnQn}).
Since $K$ is a compact set,
this infimum is reached in $K$.
If it is reached in $X \in {\cal O}$
then Lagrange's theorem gives
that $X$ is an eigenvector of $Q_n$
and $\epsilon_n$ is the associated eigenvalue.
If it is reached in $K\setminus {\cal O}$
then, by induction, $\epsilon_n$ 
is greater than or equal to
some eigenvalue of $Q_{n'}$
for some $n' < n$.
Then we just have to study 
$\Sp (Q_n)$,
the spectrum of $Q_n$,
for a generic $n$.
We claim: 
$$
\Sp (Q_n) =
\left\{ 1 - \cos\frac{(2k+1)\pi}{2n}
:\:
k \in \{0; 1; \dots; n-1\} \right\}
.$$
Indeed, with $(-\Delta_n)$
the opposite
of the discrete Laplacian
on a segment of $n$ sites
with $0$ boundary conditions,
i.e., the operator
obtained from $Q_n$
by replacing the two coefficients
$-1/\sqrt{2}$ by $-1/2$,
and with $\chi_n$
the characteristic polynomial
of $(-\Delta_n)$, 
we have for any $\lambda \in \R$,
and with $I$ the identity matrix,
$$
\det (Q_n-\lambda I)
= (1-\lambda) \chi_{n-1}(\lambda) - \frac{1}{2} \chi_{n-2}(\lambda)
$$
(set $\chi_0 := 1$,)
while, for any $k \geq 2$,
$$
\chi_k (\lambda)
= (1-\lambda) \chi_{k-1}(\lambda)
- \frac{1}{4} \chi_{k-2}(\lambda)
,$$
so that, for $0 < \lambda <2$,
$$
\chi_k(\lambda)
= \alpha \zeta^k + \bar{\alpha} \bar{\zeta}^k
$$ 
with
$$
\left\{
\begin{array}{l}
\zeta := \frac{1}{2}\left(1-\lambda + i \sqrt{1-(1-\lambda)^2}\right),\\
\alpha:= \zeta\left(\zeta - \bar{\zeta}\right)^{-1},
\end{array}
\right.
$$
and this gives,
still in the case $0 < \lambda < 2$,
and with $\theta$ in $]0,\pi[$,
defined by $e^{i\theta} = 2\zeta$,
\begin{eqnarray*}
\lambda \in \Sp (Q_n) 
&\Leftrightarrow&
(1-\lambda) \chi_{n-1} (\lambda) - \frac{1}{2} \chi_{n-2}(\lambda) =0\\
&\Leftrightarrow&
\chi_n(\lambda) = \frac{1}{4}\chi_{n-2} (\lambda)\\
&\Leftrightarrow&
\frac{1}{4}\left(e^{2i\theta}\alpha\zeta^{n-2}
+ e^{-2i\theta}\bar{\alpha}\bar{\zeta}^{n-2}\right)
= 
\frac{1}{4}\left(\alpha\zeta^{n-2}
+ \bar{\alpha}\bar{\zeta}^{n-2}\right)\\
&\Leftrightarrow&
e^{i\theta}\left(e^{i\theta}\alpha\zeta^{n-2}
- e^{-i\theta}\bar{\alpha}\bar{\zeta}^{n-2}\right)
=
e^{-i\theta}\left(e^{i\theta}\alpha\zeta^{n-2}
- e^{-i\theta}\bar{\alpha}\bar{\zeta}^{n-2}\right)\\
&\Leftrightarrow&
e^{i\theta}\alpha\zeta^{n-2} \in \R\\
&\Leftrightarrow&
e^{in\theta} \in i\R\\
&\Leftrightarrow&
\lambda = 1 - \cos\frac{(2k+1)\pi}{2n}
\quad\mbox{for some $k$ in }
\{0;1;\dots;n-1\}
.\end{eqnarray*}
In that way one gets all the eigenvalues
of $Q_n$ contained
in $]0,2[$,
but, since their number is $n$,
one gets the whole spectrum
of $Q_n$.
As a consequence
$$
\epsilon_n \geq 1-\cos\frac{\pi}{2n}
$$
(actually it is easy to prove
the equality)
and this ends the proof
of Lemma \ref{gammest}.
\qed

\section{From one to many obstacles}
\lbl{1tomany}

\subsection{Grouping the obstacles}
\lbl{grpobst}

At the end of this section
we will prove Theorems \ref{thm1} and \ref{thm2}
by induction on the number
of obstacles.
To that purpose
and before following
in the four next subsections
the four steps strategy
we described in the first section,
we introduce here
some tools
to group in a single obstacle a set of obstacles
which are ``close on a given scale $\sigma$".
Calling ${\cal R}$
the set of all finite sets
of rectangles on the plane,
we define in this subsection,
a family
$(g_\sigma)_{\sigma \geq 0}$ 
of transformations
of ${\cal R}$,
which in some sense group in single rectangles
the rectangles of an $S \in {\cal R}$
which have a distance smaller than $\sigma$
between them.
Actually these functions $g_\sigma$
are hardly more than an additional notation,
but they will be omnipresent
from this point up to the end of the work.

Given $\sigma \geq 0$ and
$$
S = \left\{R_1; R_2; \dots; R_{|S|} \right\} 
\in {\cal R}
$$
we define an equivalence relation
on $S$ as follows.
We say that two rectangles $R$ and $R'$
in $S$ are in the same equivalence class
if there exists a finite sequence
$R_1$, $R_2$, \dots, $R_k$
of rectangles in $S$
such that
$$
R=R_1,\;
R'=R_k\;
\mbox{and}\;
\forall i<k,
d_\infty(R_i,R_{i+1})<\sigma.
$$
Calling $C$ the set 
of the equivalent classes
we define 
(recall the notation of subsection \ref{stratnnot})
$$
\bar{g}_\sigma : S \in {\cal R}
\longmapsto
\left\{ \RC \left( \bigcup_{i\in c} R_i \right) \right\}
_{c \in C}
\in {\cal R} 
.$$
Since, for any $S$,
$$
|\bar{g}_\sigma (S)| \leq |S|
$$
with equality only if $\bar{g}_\sigma (S) = S$,
it is clear that the sequence of the iterates
$$
\left(\bar{g}^{(k)}_\sigma (S)\right)_{k\geq 0}
\in {\cal R}^\N
$$
is a stationary sequence, and we call $g_\sigma(S)$
its limit (for the discrete topology):
$$
g_\sigma(S)
:=
\lim_{k\rightarrow +\infty}
\bar{g}^{(k)}_\sigma (S)
.$$

We claim
(recall the notation of subsection
\ref{stratnnot})
\begin{prp}
\lbl{prpg}
For any $S$ in ${\cal R}$
and any $\sigma'\geq\sigma\geq0$, we have:
\begin{eqnarray*}
i)&&
\sum_{R\in g_\sigma(S)}
|\partial R|
\leq
\sum_{R\in S}
|\partial R|
+ 4\sigma
\left(|S|-|g_\sigma(S)|\right)
;\\
ii)&&
\hsh \bx{\underline{g_\sigma (S)}}_{\sigma'}
= \hsh \bx{\underline{S}}_{\sigma'}\\
&& \mbox{and} \quad 
\vsh \bx{\underline{g_\sigma (S)}}_{\sigma'}
= \vsh \bx{\underline{S}}_{\sigma'}
;\\
iii)&&
g_{\sigma'}\left(g_\sigma(S)\right)
= g_{\sigma'}(S)
.\end{eqnarray*}
\end{prp}

\noindent
{\bf Proof:}
Note that
$d_\infty(R_i,R_{i+1})<\sigma$
implies that we can construct
a rectangle on the plane $R''$, with side lengths
shorter than $\sigma$ and such that 
$$
\left\{
\begin{array}{l}
R_i \cup R'' \cup R_{i+1} \mbox{ is a connected set,}\\
\RC\left(R_i \cup R'' \cup R_{i+1}\right)
= \RC\left(R_i \cup R_{i+1}\right).
\end{array}
\right.
$$
From this it is easy to deduce $i)$ and $ii)$
for $\bar{g}_\sigma$, then for $g_\sigma$.

To prove $iii)$ 
observe that
$$
\underline{S}
\subset
\underline{g_\sigma(S)}
\quad\Rightarrow\quad
\underline{g_{\sigma'}(S)}
\subset
\underline{g_{\sigma'}(g_\sigma(S))}
,$$
and
$$
\underline{g_{\sigma}(S)}
\subset
\underline{g_{\sigma'}(S)}
\quad\Rightarrow\quad
\underline{g_{\sigma'}(g_\sigma(S))}
\:\subset\:
\underline{g_{\sigma'}(g_{\sigma'}(S))}
\:=\:
\underline{g_{\sigma'}(S)}
.$$
So, 
$$
\underline{g_{\sigma'}(g_\sigma(S))}
=
\underline{g_{\sigma'}(S)}
$$
and this gives, for $\sigma' >0$,
$$
g_{\sigma'}(g_\sigma(S)) =
g_{\sigma'}(S)
.$$
Since this equality is obvious in the case
$\sigma'=0$
this concludes the proof.
\qed

\subsection{Corridors and rough estimates}
\lbl{corr}
For a given finite set of rectangles on $\Z^2$
$$
\ds{S} :=
\left\{\ds{R}_1; \ds{R}_2; \dots; \ds{R}_{|\ds{S}|}
\right\}
$$
one can define or redefine the $R_i$'s,
without changing $\ds{S}$,
by
$$
R_i := \bx{\ds{R}_i}_{\frac{1}{2}}
.$$
Then, with
$$
S:= \{R_1;R_2;\dots;R_s\}
,$$
the hypothesis
of Theorem \ref{thm1}
$$
\inf_{i\neq j}
d_\infty(\ds{R}_i,\ds{R}_j)
> 3
$$
can now be written
$$
g_3(S) = S.
$$
This guarantees 
that any $\ds{R}_i$
can be bypassed without collision
by particles
using the corridor
$$
\bx{R_i}_2 \setminus \bx{R_i}_1
.$$
This is the key to the following
result (recall the definition
of the various stopping times
in subsection \ref{stratnnot}):
\begin{lmm}
\lbl{corrds}
Let $S$ be a finite set of rectangles
on the plane such that 
$$\bx{\ds{S}}_{\frac{1}{2}} = S,$$
$n$ and $p$ two integers
larger than or equal to 2 
and $\cf{z} = (z_1,\dots,z_n)$
in $(\Z^2)^n$.
\begin{itemize}
\item[i)]
If $\cf{z}$ and $S$ satisfy
the hypotheses
of Theorem \ref{thm1}
then, for any $\theta\geq 2$,
\begin{eqnarray*}
&&\ds{P}_{\cf{z}}
\left(\ds{\Tau}_c > \T{\ds{w}_S\geq \theta}\right)
\quad \geq \quad
\ds{P}_{\cf{z}}
\left(\ds{\Tau}_c > \T{\ds{w}_{g_\theta(S)}\geq \theta}\right)\\
&&\quad\geq\quad
\exp\left\{-\cst (n+p)n^2 \ln \theta\right\}
.\end{eqnarray*}
\item[ii)]
If, for some $\sigma \geq 3$, we have 
$$
\left\{
\begin{array}{l}
g_\sigma(S)=S\\
\sum_{R\in S} |\partial R| \leq p\sigma\\
|S|\leq p/4\\
\ds{w}_S(\cf{z})\geq\sigma
\end{array}
\right.
$$
then, for any $\theta\geq 2$, 
\begin{eqnarray*}
&&\ds{P}_{\cf{z}}
\left(\ds{\Tau}_c > \T{\ds{w}_S\geq \theta\sigma}\right)
\quad\geq\quad
\ds{P}_{\cf{z}}
\left(\ds{\Tau}_c > \T{\ds{w}_{g_{\theta\sigma}(S)} \geq \theta\sigma}\right)\\
&&\quad\geq\quad
\exp\left\{-\cst (n+p)n^2 \ln \theta\right\}
.\end{eqnarray*}
\end{itemize}
\end{lmm}

\noindent
{\bf Proof:}
The proof goes as follows.
We first prove {\em i)}
in the case $\theta=3$,
then adapt the proof to get {\em ii)}
in the case $\theta=3$,
we then get {\em ii)} by induction
on $\lceil \log_3 \theta \rceil$,
and finally deduce $i)$ from $ii)$
in the general case.

\medskip
{\bf First step: {\em i)} in the case $\theta= 3$.}
Assume that $\theta=3$ and that the hypotheses
of Theorem \ref{thm1} are satisfied.
In that case the first inequality in {\em i)}
is an equality, we have to prove the second inequality.
Without loss of generality
we can assume that the particles
are initially ordered
in lexicographical order
(so that $z_1$ is the most
southern of the most western
particles.) We will estimate
the probability $p_0$
of the following event
which implies
$$
\left\{
\ds{\Tau}_c > \T{\ds{w}_{g_\theta(S)}\geq \theta}
\right\}
.$$
\begin{itemize}
\item While the other particles
do not move, the first particle
moves westwards 
and uses the corridors in
$$
\bx{\underline{S}}_2
\setminus \bx{\underline{S}}_1
$$
to bypass the rectangles in $S$.
As soon as it exits from
the horizontal shadow of $\bx{g_3(S)}_3$
it stops in some $z'_1 \in \Z^2$ and we define
(recall the notation of subsection \ref{stratnnot})
$$
S'_1:=S \cup \left\{\bx{z'_1}\right\}.
$$
(Note that $g_3(S_1') = S_1'$.)
\item While the other particles
do not move, the second particle
moves westwards 
and uses the corridors in
$$
\bx{\underline{S'_1}}_2
\setminus \bx{\underline{S'_1}}_1
$$
to bypass the rectangles in $S'_1$.
As soon as it exits from
the horizontal shadow of $\bx{g_3(S'_1)}_3$
it stops in some $z'_2 \in \Z^2$ and we define
$$
S'_2:=S \cup \left\{\bx{z'_2}\right\}.
$$
\item We go on in the same way up to the last
particle's exit from the horizontal shadow
of $\bx{g_3(S'_{n-1})}_3$.
\end{itemize}
For any $k\geq 1$
the number of sites to the west of $z_k$
(on the same latitude)
contained in the horizontal shadow
of
$$
\bx{g_3(S'_{k-1})}_3
= \bx{S'_{k-1}}_3
$$
(set $S'_0 := S$)
is smaller than or equal to the total width
of this horizontal shadow,
estimated from above by 
$$
\frac{1}{2}\sum_{R\in S_{k-1}} |\partial R| + |S_{k-1}|\times 6
\leq
\frac{p+4(k-1)}{2} + 6\left(\frac{p}{4} + k-1\right)
= 2p+8(k-1)
.$$
So that, in such a scheme,
the $k^{th}$ particle
makes at most $\cst (p+k)$ steps
before stopping.
As a consequence
$$
p_0 \geq
\prod_{k=1}^{n}
\left(\frac{1}{4n}\right)
^{\cst (p+k)}
\geq
\exp\left\{-\cst (n+p)n^2\right\}
.$$

\medskip
{\bf Second step: {\em ii)} in the case $\theta= 3$.}
First of all note that the first inequality in {\em ii)}
is trivial, we just prove the second one.
If $\sigma < 4$ the previous arguments
give directly the result,
the only difference
is that we have to use
Proposition \ref{prpg}
to estimate the number of sites
to the west of a given point
and in the horizontal shadow
of some $[g_{3\sigma}(S'_{k-1})]_{3\sigma}$.
Indeed, since
$$
\hsh \bx{\underline{g_{3\sigma}(S'_{k-1})}}_{3\sigma}
= \hsh \bx{\underline{S'_{k-1}}}_{3\sigma}
,$$
this number is smaller than or equal to 
$$
\frac{1}{2}
\sum_{R\in S'_{k-1}} |\partial R|
+|S'_{k-1}|\times 6\sigma
\leq
\frac{p\sigma+4(k-1)}{2} + \left(\frac{p}{4} + k-1\right)6\sigma
\leq
\cst (p+k).
$$
If $\sigma \geq 4$
we will generalize
the previous arguments
by describing the system 
`on scale $\sigma$'.
We set
$$
\bar{\sigma} := 
\left\lfloor \frac{\sigma}{4} \right\rfloor
\geq 1
$$
and, denoting by $\ds{Z}_1$, $\ds{Z}_2$, \dots,
$\ds{Z}_n$ the processes
followed by the different particles,
we define recursively
the following stopping times
for any $k$ in $\{1;\dots;n\}$:
\begin{eqnarray*}
\Tau_{k,0} & := & 0,\\
\Tau_{k, i+1} & := &
\inf \left\{t\geq \Tau_{k,i}
:\: \bx{\ds{Z}_k(t)} \not\subset
\bx{\ds{Z}_k(\Tau_{k,i})}_{\bar{\sigma}/2} \right\}.
\end{eqnarray*}
We will say that the 
$k^{th}$ particle performs
a westward, eastward, \dots\ 
$\sigma$-step
at each time $\Tau_{k,i}$ such that the last step of
the particle was westward, eastward, \dots\
Define now 
$$
S'_0
:= S \cup \left\{\bx{z_2}_{\bar{\sigma}/2}\right\}
\cup \dots \cup
\left\{\bx{z_n}_{\bar{\sigma}/2}\right\}.
$$
and note that $g_{3\bar{\sigma}}(S'_0) = S'_0$.
One can build
a `globally westward' corridor ${\cal C}$ of width $\bar{\sigma}$,
centered at $z_1$ at its starting point,
which leads outside $\hsh \bx{g_{3\sigma}(S'_0)}_{4\sigma}$
and bypasses the rectangles in $S'_0$
using the corridors 
in 
$$\bx{\underline{S'_0}}_{2\bar{\sigma}}
\setminus \bx{\underline{S'_0}}_{\bar{\sigma}}
.$$
One can specify the orientation
of the corridor ${\cal C}$ in any of its sites
by describing ${\cal C}$ as a sequence of
westward, northward and southward
rectangular corridors,
each of them leading to the next one.
Note that if $z := \ds{Z}_1(\Tau_{1,i})$
belongs to ${\cal C}$ for some $i$,
then the probability that it performs
the next $\sigma$-step
in the direction associated to $z$
and reaching another point in ${\cal C}$
or the end of the corridor is,
by symmetry,
at least $1/8$.
If the first particle behaves
in this way at each $\sigma$-step,
using once again Proposition~\ref{prpg}
which gives that the corridor ${\cal C}$
has a `length' smaller than or equal to
$\cst (p+n)\sigma$,
we get that it follows the
the whole corridor in not more
than $\cst(p+n)$ $\sigma$-steps
and remaining confined inside
$$
\bx{{\cal C}}_{\bar{\sigma}/2}
= \bigcup_{z \in C}
\bx{z}_{\bar{\sigma}/2}
.$$
If we require also that any
$\sigma$-step of the first particle
is made in a time smaller
than $\bar{\sigma}^2$
(and, since $\bar{\sigma}^2$
is the typical order
of the time spent
to perform a $\sigma$-step,
this occurs at each time 
with a probability
which can be bounded
from below by a constant $q>0$),
then the total time spent
to follow the whole corridor,
is smaller than or equal to
$\cst\cdot (p+n)\bar{\sigma}^2$.
By Brownian approximation
and using the reflection principle,
it is easy to see that
the probability,
for any given $k\geq 2$,
that the $k^{th}$ particle
did not perform
any $\sigma$-step
in this time,
i.e., remained
confined inside $S'_0$,
is bounded from below
by $\exp\{-\cst(p+n)\}$.
Since
$$
g_{3\bar{\sigma}}(S_0') = S_0'
$$
implies that 
$$
\bx{{\cal C}}_{\frac{\bar{\sigma}}{2}}
\cap
S_0'
= \emptyset
,$$
this global event
implies that the first
particle reaches some site
$z_1'$ outside of
the horizontal shadow of
$\bx{g_{3\sigma}(S_0')}_{4\sigma}$
without any collision.
It occurs with a probability
bounded from below
by
$$
\left(
\frac{q}{8}
\right)^{\cst(p+n)}
\Big(
\exp\{-\cst(p+n)\}
\Big)^{n-1}
\geq
\exp\{-\cst (p+n)n\}
.$$
Defining
$$
S'_1
:= S \cup \left\{\bx{z'_1}_{\bar{\sigma}/2}\right\}
\cup \left\{\bx{z_3}_{\bar{\sigma}/2}\right\}
\cup \dots \cup
\left\{\bx{z_n}_{\bar{\sigma}/2}\right\}.
$$
and, as previously,
building recursively a sequence
of similar events
we get eventually
\begin{eqnlbl}{ii)a=3}
\ds{P}_{\cf{z}}
\left(\ds{\Tau}_c > \T{\ds{w}_{g_{3\sigma}(S)} \geq 3\sigma}\right)
\geq
\exp\left\{-c_1 (n+p)n^2\right\}
\end{eqnlbl}
for some constant $c_1$
independent of
$\cf{z}$, $S$, $\sigma$, $p$
and any other parameter.

\medskip
{\bf Third step: {\em ii)} in the general case.}
Define now for any $\theta\geq 2$
$$
m:= \lceil \log_3 \theta \rceil
$$ 
and
$$
S_m := g_{3^m\sigma} (S)
.$$
We will prove by induction
on $m$ that
\begin{eqnlbl}{indhyp0}
\ds{P}_{\cf{z}}
\left(\ds{\Tau}_c > \T{\ds{w}_{S_m} \geq 3^m\sigma} \right)
\geq
\exp\left\{- 2c_1 m(n+p)n^2 \right\}
.\end{eqnlbl}
We have already proved the stronger result
(\ref{ii)a=3})
in the case $m=1$,
so assume that
(\ref{indhyp0}) holds
for some $m\geq 1$.
Note that by Proposition \ref{prpg}
$$
\ds{w}_{S_{m+1}}
= \ds{w}_{g_{3.3^m\sigma}(S)}
= \ds{w}_{g_{3.3^m\sigma}(S_m)}
$$
and
$$
\left\{
\begin{array}{l}
g_{3^m\sigma}(S_m) = S_m\\
\sum_{R\in S_m} |\partial R|
\leq 2p3^m\sigma\\
|S_m| \leq p/4
\end{array}
\right.
$$
so that,
for any $\cf{z}'$ such that
$$
\ds{w}_{S_m}(\cf{z}') \geq 3^m\sigma
,$$
we have by (\ref{ii)a=3}),
applied to $\cf{z'}$, $S_m$,
$3^m\sigma$ and $2p$
instead of $\cf{z}$, $S$, $\sigma$
and $p$:
$$
\ds{P}_{\cf{z}'}
\left(\ds{\Tau}_{c;S}
> \T{\ds{w}_{S_{m+1}}\geq 3.3^m\sigma}\right)
\geq
\exp\left\{-c_1 (n+2p)n^2 \right\}
.$$
This implies, together with the strong Markov property
applied at time 
$$\T{\ds{w}_{S_m} \geq 3^m\sigma}$$
and the inductive hypothesis,
that
$$
\ds{P}_{\cf{z}}
\left(\ds{\Tau}_c 
> \T{\ds{w}_{S_{m+1}} \geq 3^{m+1}\sigma} \right)
\geq
\exp\left\{- 2c_1(m+1)(n+p)n^2 \right\}
$$
and concludes the proof of {\em ii)}.

\medskip
{\bf Fourth step: {\em i)} in the general case.}
We can get {\em i)} in the general case
as a consequence of {\em i)} in the case
$\theta=3$ and {\em ii)} in the case $\sigma = 3$
by applying the strong Markov property
at time $\T{\ds{w}_S \geq 3}$.
\qed

\medskip
Since it is straightforward
to generalize the notion of $\sigma$-step
used in the previous proof
to the continuous case of Brownian particles
it is easy to adapt this proof
to get the continuous version 
of the same results:
\begin{lmm}
\lbl{corrct}
Let $S$ be a finite set of rectangles
on the plane,
$n$ and $p$ two integers
larger than or equal to 2 
and $\cf{z} = (z_1,\dots,z_n)$
in $(\R^2)^n$.
\begin{itemize}
\item[i)] 
If $\cf{z}$ and $S$ satisfy
the hypotheses
of Theorem \ref{thm2}
then, for any $\theta\geq 2$, 
\begin{eqnarray*}
&&P_{\cf{z}}
\left(\Tau_c > \T{w_S\geq \theta}\right)
\quad\geq\quad
P_{\cf{z}}
\left(\Tau_c > \T{w_{g_\theta(S)}\geq \theta}\right)\\
&&\quad\geq\quad
\exp\left\{-\cst (n+p)n^2 \ln \theta\right\}
.\end{eqnarray*}
\item[ii)]
If, for some $\sigma \geq 3$, we have 
$$
\left\{
\begin{array}{l}
g_\sigma(S)=S\\
\sum_{R\in S} |\partial R| \leq p\sigma\\
|S|\leq p/4\\
w_S(\cf{z})\geq\sigma
\end{array}
\right.
$$
then, for any $\theta\geq 2$,
\begin{eqnarray*}
&&P_{\cf{z}}
\left(\Tau_c > \T{w_S \geq \theta\sigma}\right)
\quad\geq\quad
P_{\cf{z}}
\left(\Tau_c
> \T{w_{g_{\theta\sigma}(S)} \geq \theta\sigma}\right)\\
&&\quad\geq\quad
\exp\left\{-\cst (n+p)n^2\ln \theta\right\}
.\end{eqnarray*}
\end{itemize}
\end{lmm}

In the case of $n$ Brownian particles
and a single fixed obstacle $\pt{O}$,
it is always possible to increase
the distances between the particles
by driving away one by one the particles,
starting from the most distant from the origin
and repeating the procedure
up to the closest ones.
This allows us to release partially
the hypotheses to get a similar result.

\begin{lmm}
\lbl{corrc1obst}
For any $n\geq 2$,
any $\sigma \geq 2$
and $\cf{z}$ in $(\R^2)^n$
such that $\delta(\cf{z}) \geq \sigma$
we have, for any $\theta \geq 2$,
$$
P_{\cf{z}}
\left( \Tau_{c;\pt{O}}
> \T{\delta \geq \theta\sigma} \right)
\geq
\exp\{-\cst n^3 \ln \theta\}
.$$
\end{lmm}

The proof goes basically
in the same way as the proofs
of Lemmas \ref{corrds} and~\ref{corrct}
and we omit it.

\medskip
\noindent
{\bf Remark:}
The lower bounds appearing in this subsection
were proved by construction of suitable events.
Since, for the events we built,
the control on particles movements
is very strict
(and that is why we have only
rough estimates)
we can get as corollaries
of the proofs
some slightly stronger results.
For example we can require
not only that
$$\left\{\Tau_c > \T{w_{g_\theta(S)}\geq \theta}\right\}$$
as in Lemma \ref{corrct}-i)
but also that, for some $\alpha= \cst (n+p)$,
$$
\left\{\T{\rho\geq\alpha\theta}
> \T{w_{g_\theta(S)}\geq \theta}\right\}
,$$
without changing the lower bound we gave.
Indeed the events we built
give also a control on the maximal individual elongation,
so that, under the hypotheses of Lemma \ref{corrct}-i),
we have,
for some
$\alpha=\cst (n+p)$,
\begin{eqnlbl}{corrcti)+rho}
P_{\cf{z}}
\left(\Tau_c \wedge \T{\rho\geq\alpha\theta}
> \T{w_{g_\theta(S)}\geq \theta}\right)
\geq
\exp\left\{-\cst (n+p)n^2\ln \theta\right\}.
\end{eqnlbl}

\subsection{A logarithmic scale invariance}

In this subsection we return
to the study of the continuous
system.
Observe that
the conclusion (\ref{lnsclinv})
of our key lemma
(Lemma \ref{klmm}) shows
a logarithmic scale invariance property.
We will use this property
and the previous results
to reinforce the key lemma
by giving a lower bound
to the probability
of increasing $\delta$ from $a$
to $b$ (rather than travelling
the distance $b$ in $(\R^2)^n$)
without collision 
and without reaching
a maximal individual elongation
$\alpha b$, with the coefficient
$\alpha$ depending only on $n$.
More precisely (recall the notation
of subsection \ref{stratnnot}):
\begin{lmm}
\lbl{lsip}
There exists a positive constant
$c_0<+\infty$ such that
for any $n\geq2$,
$a \geq 2$, $\epsilon >0$
and any $\cf{z}$ in $(\R^2)^n$
such that $\delta(\cf{z}) \geq a$,
if $b\geq a^{1+\epsilon}$ then
$$
P_{\cf{z}}
\left( \Tau_{c;\pt{O}} > \T{\delta \geq b}
\mbox{ and }
\T{\rho \geq \alpha b} > \T{\delta \geq b} \right)
\geq
\left(\frac{\ln a}{\ln b}\right)^{c_\epsilon n^4}
,$$
where
$$
\alpha= c_0 n^8
$$
and $c_\epsilon$
is a constant depending only on $\epsilon$.
\end{lmm}

\noindent
{\bf Proof:}
We prove the lemma in two steps.
We first prove that
in the case $b \leq a^2$
\begin{eqnlbl}{lsipbnddb}
P_{\cf{z}}
\left( \Tau_{c;\pt{O}} 
\wedge
\T{\rho \geq \alpha b} > \T{\delta \geq b} \right)
\geq
\left(\frac{1}{2}\right)^{\cst n^4}
,\end{eqnlbl}
then we apply $\lceil \log_2 \log_a b\rceil$
times the strong Markov property to conclude.  

So, take $b \leq a^2$. For any $k >1$
we have
\begin{eqnarray*}
&&P_{\cf{z}}
\left( \Tau_{c;\pt{O}}
> \frac{b^2}{k}\right)
\quad\geq\quad
P_{\cf{z}}
\left( \Tau_{c;\pt{O}}
> \T{\partial B_2(\cf{z}, b)}
\geq \frac{b^2}{k}\right)\\
&&\quad \geq \quad
P_{\cf{z}}
\left(\Tau_{c;\pt{O}}
> \T{\partial B_2(\cf{z}, b)}\right)
-
P_{\cf{z}}
\left(\T{\partial B_2(\cf{z}, b)}
< \frac{b^2}{k} \right)
.\end{eqnarray*}
The first term
of the right hand side
can be estimated from below by
$$
P_{\cf{z}}
\left( \Tau_{c;\pt{O}}
> \T{\partial B_2(\cf{z}, a^2)}\right)
\geq
\left(\frac{1}{2}\right)^{c_1 n^4}
$$
for some constant $c_1$
given by Lemma \ref{klmm} 
with $\epsilon =1$,
while, by reflection
principle and exponential
inequality for Brownian motion,
the second one can be estimated
from above by
$$
2n \cdot 2 \exp\left\{-\frac{k}{2\cdot 2n}\right\}
\leq
\frac{1}{2}
\left(\frac{1}{2}\right)^{c_1 n^4}
$$
provided that
\begin{eqnlbl}{condt1}
k \geq \cst n^5
.\end{eqnlbl}
Now, 
\begin{eqnarray}
&&P_{\cf{z}}
\left( \Tau_{c;\pt{O}}
> \T{\delta \geq \frac{1}{k}\frac{b}{\sqrt{k}}} 
\right)
\quad\geq\quad
P_{\cf{z}}
\left(\Tau_{c;\pt{O}} 
> \frac{b^2}{k}
\geq \T{\delta\geq\frac{1}{k}\frac{b}{\sqrt{k}}}
\right)\nonumber\\
&&\quad\geq\quad
P_{\cf{z}}
\left(\Tau_{c;\pt{O}} 
> \frac{b^2}{k}\right)
-
P_{\cf{z}}
\left(\frac{b^2}{k}
< \T{\delta\geq\frac{1}{k}\frac{b}{\sqrt{k}}}
\right)
.
\label{**}
\end{eqnarray}
If (\ref{condt1}) holds
then the first term of the right hand side
(\ref{**})
can be estimated from below by
$$
\frac{1}{2}
\left(\frac{1}{2}\right)
^{c_1 n^4}
\geq
\left(\frac{1}{2}\right)
^{c_2 n^4}
$$
for some constant $c_2$,
while, dividing the time
$b^2/k$
into $k$ intervals
of length
$$
T':= \frac{b^2}{k^2}
,$$
observing that, by scaling invariance,
for any $\cf{z}'$
\begin{eqnarray*}
&& P_{\cf{z}'}
\left(\delta(\cf{Z}(T'))
< \frac{1}{k}\frac{b}{\sqrt{k}}\right)
\quad = \quad
P_{\cf{z}'}
\left(\delta(\cf{Z}(T'))
< \frac{1}{\sqrt{k}}\sqrt{T'}\right)\\
&&\quad\leq\quad
\left(n+\frac{n(n-1)}{2}\right)
\left(\frac{1}{\sqrt{2\pi}}
\frac{2}{\sqrt{k}}\right)^2
\quad\leq\quad
\frac{n^2}{k}
,\end{eqnarray*}
and using the Markov property,
the second term
of the right hand side (\ref{**})
can be estimated
from above by
$$
\left(\frac{n^2}{k}\right)
^k
\leq
\frac{1}{2}
\left(\frac{1}{2}\right)^{c_2 n^4}
$$
provided that
\begin{eqnlbl}{condt2}
k\geq \cst n^4
.\end{eqnlbl}
Choose $k = \cst n^5$ in order
to have (\ref{condt1})
and (\ref{condt2}) satisfied.
If
$$
\sigma:=\frac{b}{k\sqrt{k}}\geq 2
$$
then, applying the strong Markov property
at time $\T{\delta\geq \sigma}$
and Lemma \ref{corrc1obst}
with $\theta:= k^{3/2}$,
we get
$$
P_{\cf{z}}
\left( \Tau_{c;\pt{O}}
\geq \T{\delta\geq b}
\right)
\geq 
\frac{1}{2}
\left(\frac{1}{2}\right)^{c_2 n^4}
\exp\left\{-\cst n^2 \ln n \ln k^{3/2}\right\}
$$
so that
\begin{eqnlbl}{atob}
P_{\cf{z}}
\left( \Tau_{c;\pt{O}}
\geq \T{\delta\geq b}
\right)
\geq
\left(\frac{1}{2}\right)^{c_3 n^4}
\end{eqnlbl}
for some constant $c_3$.
If, on the contrary,
$$
b < 2k^{3/2}
$$
then (\ref{atob}) is a direct
consequence of Lemma \ref{corrc1obst}
applied to
$$
\sigma:= a\geq 2
$$
and
$$
\theta:=\frac{b}{a} \leq k^{3/2}
.$$
Finally, for any $\alpha>1$, 
\begin{eqnarray*}
&&P_{\cf{z}}
\left(\Tau_{c;\pt{O}}
> \T{\delta \geq b},\:
\T{\rho\geq\alpha b}
> \T{\delta \geq b}\right)\\
&&\quad\geq\quad
P_{\cf{z}}
\left(\Tau_{c;\pt{O}}
> \T{\delta \geq b},\:
\alpha b^2
> \T{\delta \geq b},\:
\T{\rho\geq\alpha b}
>\alpha b^2\right)\\
&&\quad\geq\quad
P_{\cf{z}}
\left(\Tau_{c;\pt{O}}
> \T{\delta \geq b}\right)
-
P_{\cf{z}}
\left(\alpha b^2
\leq \T{\delta \geq b}\right)
-
P_{\cf{z}}
\left(\T{\rho\geq\alpha b}
\leq\alpha b^2\right)
\end{eqnarray*}
and, like previously,
the second term of the right hand side
can be estimated from above by
$$
\left(\frac{n^2}{\sqrt{\alpha}}\right)^{\sqrt{\alpha}}
\leq \frac{1}{4}
\left(\frac{1}{2}\right)^{c_3 n^4}
$$
provided that
\begin{eqnlbl}{condt3}
\alpha \geq \cst n^8
,\end{eqnlbl}
the last term can be estimated
from above (by reflection principle
and exponential inequality) by
$$
2n e^{-\alpha/4}
\leq \frac{1}{4}
\left(\frac{1}{2}\right)^{c_3 n^4}
$$
provided that
\begin{eqnlbl}{condt4}
\alpha \geq \cst n^4
,\end{eqnlbl}
and this, with (\ref{atob}),
gives (\ref{lsipbnddb}), provided
that $\alpha = \cst n^8$ is such that
(\ref{condt3}) and (\ref{condt4})
hold.

To prove the result in the general case
we apply the strong Markov
property at times
$$
\Tau_0 := 0,\;
\Tau_1 := \T{\delta\geq a^2},\;
\Tau_2 := \T{\delta\geq a^4},\;
\dots,\;
\Tau_{m-1} := \T{\delta \geq a^{2^{m-1}}}
$$
and (\ref{lsipbnddb}) with $(a,b)$
replaced by 
$$
(a,a^2),\;
(a^2, a^4),\;
\dots,\;
(a^{2^{m-1}},b),
$$
where
$$
m:=\left\lceil
\log_2 \log_a b
\right\rceil.
$$
If in each interval
$[\Tau_i, \Tau_{i+1}]$
(set $\Tau_m := \T{\delta\geq b}$)
the maximal individual elongation
is smaller than
$\alpha (a^{2^{i+1}}\wedge b)$
then, on the whole interval
$\left[0,\T{\delta\geq b}\right]$
the maximal individual elongation
is bounded from above by
$$
\alpha a^2 + \alpha a^4 + \cdots
+ \alpha a^{2^{m-1}}+ \alpha b
\leq 2\alpha a^{2^{m-1}} + \alpha b
\leq 3\alpha b.
$$
So that we get
$$
P_{\cf{z}}
\left( \Tau_{c;\pt{O}} > \T{\delta \geq b}
\mbox{ and }
\T{\rho \geq 3\alpha b} > \T{\delta \geq b} \right)
\geq
\left(\frac{1}{2}\right)^{\cst n^4 m}
.$$
Under the hypothesis $b\geq a^{1+\epsilon}$,
the right hand side can be
estimated from below
by 
$$
\left(
\frac{\ln a}{\ln b}
\right)
^{c_\epsilon n^4}
$$
where $c_\epsilon$
is a constant depending $\epsilon$ only,
and this concludes the proof.
\qed 

\subsection{Proof of Theorem \ref{thm2}}

We prove now Theorem \ref{thm2}.
So, we take $S$ in ${\cal R}$ and $\cf{z}^0$
in $(\R^2)^n$
satisfying the hypotheses
of the theorem,
and,
as first step,
we prove by induction on 
$$
s:= |S|
$$
that for some constant
$c_1$ that we will specify later
\begin{eqnlbl}{indhyp}
\forall T\geq 4,
\quad
P_{\cf{z}^0}
\left(
\Tau_c > \T{w \geq T-1}
\right)
\geq
\left(
\frac{1}{\ln T}
\right)
^{c_1 (s+1) (p\ln p) n^4}
.\end{eqnlbl}

Clearly we just have
to deal with the case
$s=1$ to prove
(\ref{indhyp})
for both cases
$s=0$ and $s=1$.
So, assume that $S$
is a single rectangle
on the plane
with perimeter
$$|\partial S| \leq p.$$
By Lemma \ref{corrct}
and using the remark
at the end of section \ref{corr}
we have (see (\ref{corrcti)+rho})):
\begin{eqnlbl}{wgeqp}
\begin{array}{l}
P_{\cf{z}^0}
\Big(\Tau_c \wedge \T{\rho\geq\cst (n+p)p}
> \T{w_S \geq p}\Big)\\
\quad\geq\quad
\exp\Big\{-\cst (n+p)n^2 \ln p\Big\}
.
\end{array}
\end{eqnlbl}
We take now $O$, the origin of the plane,
at the center of the rectangle,
define
$$
B:=\overline{B_2\left(O,p/4\right)},
$$
and observe that
$S \subset B$
and that for any $\cf{z}$ in $(\R^2)^n$
$$
w_S\left(\cf{z}\right) \geq p
\;\Rightarrow\;
\delta\left(\cf{z}\right) \geq p.
$$
Consider now $\Tau_{c;B}^{p/2}$,
the first collision time
for one fixed particle $B$
and $n$ Brownian particles
with the same diameter
$$
{\rm diam}\,B = \frac{p}{2}
$$
and starting in a configuration
$\cf{z}$ in $(\R^2)^n$
such that
$$
\delta(\cf{z}) \geq 2 \frac{p}{2}.
$$
We will certainly have
$\Tau_c\geq\Tau_{c;B}^{p/2}$
and a homothety of coefficient
$2/p$ allows us to use Lemma \ref{lsip}
with $a=2$ and $b=T\geq a^2$
to get 
$$
P_{\cf{z}}
\left(\Tau_{c;S} > \T{\delta \geq \frac{p}{2}T}
\mbox{ and }
\T{\rho \geq \alpha \frac{p}{2} T}
> \T{\delta \geq \frac{p}{2}T}
\right)
\geq
\left(\frac{1}{\ln T}\right)
^{\cst n^4}
.$$
Observing that for any $\cf{z}$
$$ \delta (\cf{z}) \geq \frac{p}{2}T
\;\Rightarrow\;
w_S(\cf{z}) \geq T-1,$$
and using the strong Markov
property at time
$\T{w_S\geq p}$ to combine
this last result
with (\ref{wgeqp}),
we conclude:
\begin{lmm}
\lbl{sihs=1}
There exist
an $\alpha= \cst (p+n^8)>1$
and a constant $c_2<+\infty$,
such that for any
single rectangle on the plane $S$
and any $\cf{z}$ satisfying
the hypotheses of Theorem~\ref{thm2},
we have, for any $T\geq 4$,
$$
P_{\cf{z}}
\left(
\Tau_{c;S} \wedge \T{\rho\geq \alpha\frac{p}{2}T} 
>\T{w_S\geq T-1}
\right)
\geq
\left(\frac{1}{\ln T}\right)
^{c_2 (p\ln p) n^4}
.$$
\end{lmm}
This result implies
(\ref{indhyp}) for
$s=0$ and $s=1$,
provided that
\begin{eqnlbl}{c1condt1}
c_1\geq c_2/2.
\end{eqnlbl}
This is the first condition
to determine the choice
of $c_1$
and we will assume in the sequel
that it is satisfied.

For $s\geq 2$, we prove now
(\ref{indhyp}) assuming
the result
for any smaller $s$.
To that purpose, given $T\geq 4$,
we distinguish two cases.

\noindent
{\bf Case 1:} None of the $k$ connected components
of 
$$
D:=\bigcup_{i=1}^n
B_2 \left(
z^0_i, \alpha p T\right)
,$$
where $\alpha$ is like in the previous Lemma,
intersects more than one rectangle of $S$.
In that case, applying Lemma \ref{sihs=1}
to the $k$ systems formed by the
$n_j$ particles contained in the $j^{th}$
connected component of $D$ (with $1\leq j \leq k$)
we get
$$
P_{\cf{z}^0}
\left(
\Tau_c > \T{w \geq T-1}
\right)
\geq
\prod_{j=1}^k
\exp
\left\{
- c_2 (p\ln p) n_j^4 \ln \ln T
\right\}
$$
and, since
$$
n_1 + n_2 + \cdots + n_k = n
,$$
this gives (\ref{indhyp}). 

\noindent
{\bf Case 2:}
One of the connected components
of $D$ intersects more than
one rectangle of $S$.
In that case we introduce a `mesoscopic scale'
$$
\sigma_0 :=
\inf\left\{
\sigma\geq 3 :\:
|g_\sigma(S)|<|S|\right\}
$$
which lies
between the `microscopic scale'
1 which is the diameter of the particles,
and the `macroscopic scale' ${2pn\alpha}T$
(of order $T$ for large $T$)
as a consequence of our case 2 hypothesis.
Assume that 
$$
T':= \frac{1}{2pn\alpha}\sigma_0
$$
is larger than or equal to 4
(we will soon explain 
why this hypothesis is not restrictive),
then,
like in case 1,
Lemma \ref{sihs=1}
applied to $T'$ instead of $T$ 
gives
$$
P_{\cf{z}^0}
\left(
\Tau_c > \T{w_S \geq T' -1}
\right)
\geq
\exp\left\{-c_2 (p\ln p) n^4 \ln\ln T\right\}
.$$
Now if $\cf{z}^1$ in $(\R^2)^n$
is such that $w_S(\cf{z}^1)\geq T'- 1$,
then by Lemma \ref{corrct}-ii) applied
to $T' -1$ instead of $\sigma$
(note that, by construction, $g_{T'-1}(S) = S$)
and $\theta := 4pn\alpha$ gives
$$
P_{\cf{z}^1}
\left(\Tau_c
> \T{w_{g_{\sigma_0}(S)} \geq \sigma_0}\right)
\geq
\exp\left\{-\cst (n+p)n^2 \ln (pn)\right\}.
$$
Combining these last two estimates
with the strong Markov property
at time $\T{w_S\geq T'-1}$
we get that there is a constant
$c_3$ independent of any parameter,
such that
$$
P_{\cf{z}^0}
\left(\Tau_c
> \T{w_{g_{\sigma_0}(S)} \geq \sigma_0}\right)
\geq
\exp\left\{-c_3 (p\ln p) n^4 \ln\ln T\right\}
.$$
And the constant $c_3$ can be chosen
to cover also the case $T'<4$:
this is a consequence of Lemma \ref{corrct}-i).

Define now,
for any $k\geq 0$,
$$
S_k:= g_{4^k \sigma_0}(S)
,$$
define
$$
k_0:= \inf\left\{k\geq 0 :\: S_k= S_{k+1}\right\} 
,$$
and set
$$
\sigma_1 := 4^{k_0}\sigma_0 
.$$
It is easy to check
$$
|g_{\sigma_1}(S)| \leq s-1-k_0
$$
and, using once again Lemma \ref{corrct}
(with $\sigma = \sigma_0$ and $\theta = 4^{k_0 +1}$)
and the previous estimate,
we get 
$$
P_{\cf{z}^0}
\left(\Tau_c
> \T{w_{g_{\sigma_1}(S)} \geq 4 \sigma_1}\right)
\geq
\exp\left\{-c_4 (k_0+1) (p\ln p) n^4 \ln\ln T\right\}
$$
for some constant $c_4$ independent
of any parameter.
Considering, like previously,
the first collision time
for larger particles
of diameter $\sigma_1$,
initially centered at some $\cf{z}$
such that
$w_{g_{\sigma_1}(S)} \geq 4 \sigma_1$
and with $g_{\sigma_1}(S)$
as set of fixed obstacles,
using the strong Markov property at time
$$\T{w_{g_{\sigma_1}(S)} \geq 4 \sigma_1},$$
an homothety of coefficient $1/\sigma_1$,
and the inductive hypothesis,
we conclude
\begin{eqnarray*}
&&P_{\cf{z}^0}
\left(\Tau_c
> \T{w_S \geq T-1}\right)\\
&&\quad\geq\quad
\exp\left\{-c_4 (k_0+1) (p\ln p) n^4 \ln\ln T
-c_1 (s-k_0) (p\ln p) n^4 \ln\ln T\right\}
.\end{eqnarray*}
This implies (\ref{indhyp})
provided
$$
c_1 \geq c_4
$$
which, after (\ref{c1condt1}),
is our only constraint on $c_1$,
and this concludes our first step.

As a consequence of this result,
since $|S|\leq p/4$ we have
$$
\forall T\geq 4,
\quad
P_{\cf{z}^0}
\left(
\Tau_c > \T{w \geq T-1}
\right)
\geq
\left(
\frac{1}{\ln T}
\right)
^{\cst (p^2\ln p) n^4}
.$$
So, to conclude the proof
of Theorem \ref{thm2}
we just have to notice
that if $w$ has been increased
without collision
up to $T-1$,
then there cannot be any collision
before time $T$
unless some particles
have a superdiffusive behaviour.
By the exponential inequality, for $T\geq 4$:
\begin{eqnarray*}
P_{\cf{z}^0}\left(
\Tau_c > T\right)
& \geq & 
P_{\cf{z}^0}\left(
\Tau_c > \T{w \geq T-1},
\Tau_c > T\right)\\
&\geq&
\left(\frac{1}{\ln T}\right)
^{\cst(p^2\ln p)n^4}
- n^2 2\exp\left\{
-\frac{((T-1)/\sqrt{2})^2}{2(2T)}
\right\}
.\end{eqnarray*}
This last expression
can be estimated
from below by
$$
\left(\frac{1}{\ln T}\right)
^{\cst(p^2\ln p)n^4}
$$
provided that
$$
T\geq (\cst (p^2\ln p) n^4)^2,
$$
and this ends the proof.
\qed

\subsection{Proof of Theorem \ref{thm1}}

We deduce Theorem \ref{thm1} from
Theorem \ref{thm2}
using the strong coupling
established
by Koml\'os, Major and Tusn\'ady:

\begin{thm}[Koml\'os - Major - Tusn\'ady]
There exist three constants
$\lambda'$, $C$, and $K$ in $]0,+\infty[$
and there exists a probability space
$(\Omega, {\cal F}, P)$
on which can be defined, for any $n\geq 1$,
$n$ two-dimensional independent Brownian motions
$Z_1$, $Z_2$,~\dots, $Z_n$
and $n$ two-dimensional independent
continuous time random walks
$\ds{Z}_1$, $\ds{Z}_2$,~\dots, $\ds{Z}_n$
with $\cf{Z}(0) = \ds{\cf{Z}}(0)$,
such that for all $x > 0$
and each $T \geq 0$
$$
P\left(
\sup_{1\leq k \leq n}
\sup_{0\leq t \leq T}
\left\|\frac{1}{\sqrt{2}}Z_k(t) - \ds{Z}_k(t)\right\|_2
> C \ln T + x
\right)
< n K e^{-\lambda' x}
.$$
In particular there exist two constants
$C$ and $\lambda$
such that
$$
\forall T \geq 2, \quad
P\left(
\sup_{1\leq k \leq n}
\sup_{0\leq t \leq T}
\left\|\frac{1}{\sqrt{2}}Z_k(t) - \ds{Z}_k(t)\right\|_2
> (1+C) \ln T
\right)
< n e^{-\lambda \ln T}
.$$
\end{thm}

The proof of the one-dimensional version,
in the case $n=1$,
of this theorem is given in
\cite{KMT1} and \cite{KMT2},
and the generalization
to the two-dimensional situation
with $n\geq 1$ is straightforward.
This theorem implies that
with high probability
the particles performing
random walks remain
contained up to time $T\geq 4$ in
balls of diameter
$$\sigma_0:=3(1+C)\ln T,$$
centered at some 
rescaled Brownian
motions  $Z_k/\sqrt{2}$.

A way to realize the event
$\{\ds{\Tau}_{c;\ds{S}}>T\}$
is to reach,
without collision,
a configuration $\cf{z}^1$
such that the particles
are at distance $5\sigma_0$,
at least,
one from each other and each from
some $g_{5^k\sigma_0}(S)$
(once again we assume $\bx{\ds{S}}_{1/2}= S$)
satisfying
\begin{eqnlbl}{kcondt}
g_{5^{k+1}\sigma_0}(S) 
= g_{5^{k}\sigma_0}(S) 
,\end{eqnlbl}
then to require
that up to time $T$
there is no collision
neither
between the balls centered at
the rescaled
Brownian motions $\cf{Z}/\sqrt{2}$
coupled with $\ds{\cf{Z}}$ and
initially starting in $\cf{z}^1$,
nor between these balls
and $g_{5^k\sigma_0}(S)$.
The probability
of the last part
of this event
can be estimated,
after homothety
of coefficient
$1/\sigma_0$,
by Theorem~\ref{thm2},
and,
since condition (\ref{kcondt})
is clearly satisfied
by some
$$
k\leq |S|\leq \frac{p}{4}
\quad,$$
the probability of the first part,
i.e.,
to reach such a configuration 
$\cf{z}^1$ without collision,
can be estimated from below
using Lemma \ref{corrds}.
So that, using the strong coupling,
we get, for $T\geq T_0 = \nu^2$
as in Theorem \ref{thm2}:
$$
\begin{array}{rccl}
\ds{P}_{\cf{z}}\left(\ds{\Tau}_c>T\right)
& \geq &&
\exp\left\{
-\cst (n+p)n^2
\ln \left(5^{p/4}3(1+C)\ln T\right)
\right\}\\
&& \times &
\exp\left\{-\cst (p^2\ln p)n^4 \ln\ln T \right\}\\
&& - &
n \exp\left\{-\lambda \ln T\right\}.
\end{array}
$$
And this last expression
can be estimated
from below by
$$
\exp\left\{-\cst (p^2\ln p)n^4 \ln\ln T \right\}
$$
provided that
$$
T\geq\exp\left\{\left(
\cst (p^2\ln p)n^4
\right)^2 \right\}
.$$
\qed

\section{Concluding remarks}

\subsection{Higher dimension}

As we wrote above
the behaviour
of the non-collision probability
is well known
for the one-dimensional case,
and we derived in this work
some estimates
for the two-dimensional case.
What about the higher dimensions?

For random walks in $\Z^d$
with $d\geq 3$
Wiener's test 
(see for example
Theorem~2.2.5 in \cite{L})
applied to the subset
$A$ of $(\Z^d)^n$
corresponding,
like in the proof of our key Lemma,
to collisions between particles
or particles and fixed obstacles
shows that $A$ is transient.
The method we followed
in this paper to estimate
the non-collision
probability up to time $T$,
gives in dimension $d\geq 3$
a somewhat stronger result,
at least in the continuous version
of the problem:
we obtain,
for the system starting in $\cf{z}$,
a lower bound
{\em depending only on $w(\cf{z})$}
(defined like in the two-dimensional case) 
to the non-collision probability
up to time $T=+\infty$.
For example, in absence
of fixed obstacles
and in the case
of Brownian particles starting
from a configuration $\cf{z}$
such that the centers of the particles
are at least distant of $a\geq 1$,
we get
(following the proof of the key lemma) 
$$
P_{\cf{z}} \left(
\Tau_c=+\infty \right)
\geq
\left(1-\frac{1}{a^{d-2}}\right)
^{\frac{{\scriptstyle n(n-1)}}
{{\scriptstyle 2\left(1-
\cos\frac{{\scriptstyle \pi}}
{{\scriptstyle n+1}}\right)
}}
}
\quad.$$
\noindent
{\bf Remark:}
The reason why we get this $\cos(\pi/(n+1))$
instead of $\cos(\pi/2n)$
like in Lemma \ref{klmm}
which deals with the case
of Brownian particles
{\em with one fixed obstacle},
is that in absence of fixed obstacles
we have to study the spectrum
of the opposite
of the discrete Laplacian
$(-\Delta_n)$
instead of the spectrum
of the operator $Q_n$
we introduced
in the proof of the key lemma.
This is not specific of the dimension $d$,
things go in the same way
in dimension 2 when there are no fixed obstacles.

\subsection{What is the right exponent?}

We have proven a lower bound
to the non-collision probability
up to time $T$.
The question we will address
in this last subsection
is the question of the `right'
asymptotics for large $T$.
To make the problem simpler,
let us first consider
the non-collision probability
without fixed obstacles,
i.e., the case $p=2$.
Since the difference of two random walks
(or two Brownian motions)
is a rescaled random walk
(or Brownian motion)
in the case $n=2$
we have
$$
P_{\cf{z}}(\Tau_c>T)
\sim \frac{\cst\!(\cf{z})}{\ln T}
\quad,
$$
with $\cst\!(\cf{z})$ a constant
which depends on $\cf{z}$,
and since, for $n>2$,
the collisions between the first and second particle,
the third and fourth particle, and so on
are independent events,
we have,
for any $T$ larger than some $T_0 (\cf{z})$,
$$
P_{\cf{z}}\left(\Tau_c>T\right)
\leq {\cst\!(\cf{z})}
\left(\frac{1}{\ln T}\right)
^{\left\lfloor\frac{n}{2}\right\rfloor}
.$$  
So that the `right' asymptotic
lies between two powers of the inverse
of $\ln T$, one of which goes like $n$
and the other one like $n^4$.

This $n^4$ comes from the estimates
of our key lemma,
and the accuracy of these estimates
has to be discussed along two fault lines.
The first one is the global method
we followed: we looked for some
subharmonic function to estimate
an harmonic one linked to
our non-collision probability,
and one can discuss the form
under which we looked for this
subharmonic function.
The second fault line is the fact that
we made a very rough estimate
using $\gamma -1$ 
as lower bound
to the collision correlation.
The other estimates we made
are quite precise and most of the inequalities
we wrote are actually equalities.
This led us for some time
during the redaction of this paper
to begin to think that this behaviour
in $n^4$ was not so far
from the `right' estimate,
and to doubt about the accuracy
of the power $n(n-1)/2$ 
that one expects
(for example by analogy
with the one-dimensional case.) 
But we performed some numerical simulations
which tend to show
that the non-collision probability
up to $T$
behaves like
$$
\frac{\cst\!(\cf{z})}
{(\ln T)^{\nu(T)}}
$$
where $\nu(T)$
is a function
which {\em grows}
slowly towards $n(n-1)/2$.

Note that even if $n(n-1)/2$
can be imagined as the consequence
of some decorrelation for large $T$
between the collision regarding
the different pairs of particles,
it is easy to see, that,
at least in the case of dimension 1,
there is no such decorrelation: 
up to the first collision the particles
keep their initial ordering
and there is no decorrelation between,
say, the collisions regarding the first
three particles.

Using this observation on the conserved
ordering in dimension one
(which implies that the possible collisions
are $(n-1)$ and not anymore $n(n-1)/2$),
the method we followed in this paper,
would have give,
for this one dimensional case,
a power (of $1/\sqrt{T}$ and
not anymore $1/\ln T$) which goes like $n^3$,
i.e., one order higher than the correct answer.
Like written above, the correct exponent $n(n-1)/2$
is given in dimension 1 by a reflection argument,
that cannot be extend, at least directly,
to higher dimensions.
But the same reflection argument
can be used, as in \cite{KOR}, to prove that
$$
(x_1, x_2, \dots, x_n) \in \R^n
\mapsto
\prod_{i<j} (x_j - x_i)
$$
is harmonic,
and it is easy to get the right exponent
from this result.
In dimension~2 the corresponding function
would be
$$
h:
(z_1, z_2, \dots, z_n) \in (\R^2)^n
\mapsto
\prod_{i<j} \ln\left\|z_j - z_i\right\|_2
,
$$
i.e., with the notation 
of the proof of the key Lemma,
$$
h:= \prod_k \ln \left(\sqrt{2} r_k\right)
= \prod_k \ln \left(\alpha_k r_k\right)
,$$
with $k$ describing
the set of indices 
such that the associated
subspace $F_k$ is of the `first kind'.
Would $h$ be harmonic or subharmonic
where it is positive,
we would get the exponent $n(n-1)/2$
in the same way we get
the exponent going like $n^4$
in the proof  of the key lemma,
using this function $h$
instead of the function $g$
we built.
Unfortunately $h$ {\em is not}
subharmonic.
But it might be possible,
to improve our result
using similar ideas
with quite precise estimates
of $\Delta h$.
Another way of improving
our result
could be based on the construction
of a subharmonic function $g$
of the form
$$
g:=\prod_k \left(\ln \alpha_k r_k\right)^{1/\gamma}
$$
with $\gamma$ a non-trivial function
(in the proof of the key Lemma
we built such a $g$
with $\gamma$ a constant depending on $n$.)

As far as the question
of `the right exponent in $p$'
is concerned,
we think that our lower bound
could be improved
up to obtainment of an exponent
independent of $p$, i.e.,
an estimate of the kind:
\begin{eqnlbl}{prot}
\forall T\geq T_0(n,p), \quad
P_{\cf{z}}\left(\Tau_c\geq T\right)
\geq
\left(\frac{c(p)}{\ln T}\right)
^{\nu(n)}
.\end{eqnlbl}
Indeed in the simpler case
of a single Brownian particle
evolving between $s$ fixed
particles,
writing $r_k(z)$
the distance between $z$ in $\R^2$
and the center
of the $k^{th}$ fixed particle,
defining the harmonic function
$$
h:= \sum_k \ln r_k
$$
and calling $\lambda$
the supremum of $h$
on $\underline{S}$,
the part of the plane
occupied by the fixed particles,
we have that $\underline{S}$
is contained in
$$
A:= \left\{z \in \R^2 :\:
h(z) \leq \lambda\right\}
$$
and, since $h$ is harmonic,
it is easy to estimate,
for $z$ such that $h(z) > \lambda$,
the probability
$P_z(\T{A} > T)$.
Observing that
the more the fixed particles
are distant from each other,
the more $A$ fits $\underline{S}$,
it is then easy
to get, in that case,
an estimate like (\ref{prot}).

Our estimates can then certainly
be improved.
But our original motivations
(see subsection \ref{mtvs})
just required estimates
going like the inverse
of `some' power of $\ln T$
for the discrete non-collision
probability.
And this is what gives
Theorem~\ref{thm1}.
  
\subsection*{Acknowledgments}
I thank Wendelin Werner
for having revealed to me
the role the logarithmic scale invariance
could play in the extension
of the result from one to many obstacles,
and for having presented to me
the strong coupling between random walk
and Brownian motion.
I thank Gabriella Tarantello
for having introduced me  
to some basic tricks of harmonic analysis,
and for having transferred to me
the conviction 
that I could drastically improve
some earlier estimates I had derived:
she was right.
I thank Salvatore Pontarelli
for having performed
the numerical simulations
which cleared up many of my doubts
on what could be `the right exponent'.
I thank Frank den Hollander
for having given to me
the opportunity to talk 
about these arguments
at Eurandom:
this was the source
of many stimulating 
and encouraging discussions.
I thank Wolfgang K\"onig
for having explained to me
what was known, unknown and believed
about the non-collision probability.
I thank Rapha\"el Cerf
for his support
and the trust he showed to me.
I thank Francesca Nardi,
Benedetto Scoppola,
Koli Ndreca and Gianluca Guadagni
for their continuous
availability
to answer or try to answer
any question I could ask them
to push forward the research
made for this paper.
I thank Enzo Olivieri and Betta Scoppola
for their continuous support
during this research,
especially in the long months
during which I entered
so many dead-ends to get these estimates.
The research in this paper
was partially supported
by Cofinanziamento 2004
Sistemi Dinamici, Meccanica Statistica
e Teoria dei Campi.

\end{document}